\documentclass{autart}
\usepackage{natbib}

\usepackage{amsmath,amssymb,theorem}
\usepackage{graphicx}
\usepackage{algorithm,algorithmic}
\usepackage{psfrag}
\usepackage{url}
\usepackage{csquotes}

\usepackage{color,xspace}
\usepackage{xcolor}
\usepackage{subfigure}
\newtheorem{theorem}{Theorem}[section]
\newtheorem{lemma}[theorem]{Lemma}

\newtheorem{remark}[theorem]{Remark}
	
\newtheorem{problem}[theorem]{Problem}

\newtheorem{definition}[theorem]{Definition}

\newcommand{\real}{{\mathbb{R}}}
\newcommand{\realpositive}{\mathbb{R}_{>0}}
\newcommand{\realnonnegative}{\mathbb{R}_{\ge 0}}
\newcommand{\integernonnegative}{\mathbb{Z}_{\ge 0}}
\newcommand{\integerpositive}{\mathbb{Z}_{> 0}}

\newcommand{\GG}{{\mathcal{G}}}

\newcommand{\NN}{{\mathcal{N}}}
\newcommand{\BB}{{\mathcal{B}}}
\newcommand{\douti}{d_i^{\operatorname{out}}}
\newcommand{\Nouti}{\mathcal{N}_i^{\operatorname{out}}}

\newcommand{\Nini}{\mathcal{N}_i^{\operatorname{in}}}
\newcommand{\Dout}{D^{\operatorname{out}}}
\newcommand{\Din}{D^{\operatorname{in}}}
\newcommand{\LL}{{\mathcal{L}}}

\newcommand{\ones}[1]{\mathbf{1}_{#1}}
\newcommand{\zeros}[1]{\mathbf{0}_{#1}}

\newcommand{\diag}[1]{\operatorname{diag}\left( #1\right)}

\newcommand{\vertices}{V}

\renewcommand{\epsilon}{\varepsilon}

\newcommand{\until}[1]{\{1,\dots, #1\}}
\newcommand{\Norm}[1]{\|#1\|}
\newcommand{\Sym}[1]{{#1}_{s}}

\newcommand{\setdef}[2]{\{#1 \; | \; #2\}}

\newcommand{\commgraph}{\GG_\text{comm}}

\newcommand{\map}[3]{#1: #2 \rightarrow #3}

\newcommand{\TwoNorm}[1]{\|#1\|}

\renewcommand{\hat}{\widehat}

\newcommand{\oprocendsymbol}{\hbox{$\bullet$}}
\newcommand{\oprocend}{\relax\ifmmode\else\unskip\hfill\fi\oprocendsymbol}

\newcommand{\longthmtitle}[1]{\mbox{}\textup{\textbf{(#1)}}}

\renewcommand{\theenumi}{(\roman{enumi})}

\newcommand{\changes}[1]{{\color{blue} #1}}
 \renewcommand{\changes}[1]{#1}

\graphicspath{{figs/}}

\begin{document}

\clearpage

\runauthor{C. Nowzari, E. Garcia, J. Cort\'es}

\begin{frontmatter}
  
  \title{Event-Triggered Communication and Control of
    \\
    {Networked} Systems for Multi-Agent Consensus}
  
  \author[address]{Cameron Nowzari} \ead{cnowzari@gmu.edu}\qquad
  \author[address2]{Eloy Garcia} \ead{eloy.garcia.2@us.af.mil} \qquad
  \author[address3]{Jorge Cort\'es} \ead{cortes@ucsd.edu}
  
  \address[address]{Department of Electrical and Computer Engineering,
    George Mason University, Fairfax, VA, 22030, USA}
  \address[address2]{Control Science Center of Excellence, Air Force
    Research Laboratory, Wright-Patterson AFB, OH, 45433, USA}
  \address[address3]{Department of Mechanical and Aerospace
    Engineering, University of California, San Diego, CA, 92093, USA}

  \begin{abstract}
    This article provides an introduction to event-triggered coordination for
    multi-agent average consensus. We provide a comprehensive account
    of the motivations behind the use of event-triggered strategies
    for consensus, the methods for algorithm synthesis, the technical
    challenges involved in establishing desirable properties of the
    resulting implementations, and their applications in distributed
    control.  We pay special attention to the assumptions on the
    capabilities of the network agents and the resulting features of the
    algorithm execution, including the interconnection topology, the
    evaluation of triggers, and the role of imperfect information.
    The issues raised in our discussion transcend the specific
    consensus problem and are indeed characteristic of cooperative
    algorithms for networked systems that solve other
    coordination tasks. As our discussion progresses, we make these
    connections clear, highlighting general challenges and tools to
    address them widespread in the event-triggered control of {networked}
    systems.
  \end{abstract}

  
\end{frontmatter}

\section{Introduction} \label{sec:intro}

This article provides an introduction to the topic of event-triggered
coordination of networked systems, with a particular emphasis on the
multi-agent consensus problem. In many applications involving multiple
agents such as vehicles, sensors, computers, etc., a group of agents
must agree upon various physical and virtual quantities of interest,
at which point it is said that a consensus has been
achieved. Consensus is a long-standing area of research, particularly
in computer science, and arguably forms the foundation for distributed
computing in
general~\citep{ROS-RMM:03c,WR-RWB:08,NAL:97,WR-RWB-EMA:07}. Consequently,
there is a vast amount of literature available on consensus
problems. \cite{ROS-JAF-RMM:07} provided a brief history of the field
and its success 10 years ago which, in a nutshell, has to do with the
extremely wide applicability of such problems across many different
disciplines. As a result, consensus problems in general are indeed
still an extremely active area of research.

On the other hand, the idea of event-triggered control has an
interesting history that only recently seems to be gaining popularity
throughout the controls community.  The basic idea of event-triggered
control is to abandon the paradigm of periodic (or continuous)
sampling/control in exchange for deliberate, opportunistic aperiodic
sampling/control to improve efficiency. For instance, it may not make
sense to constantly monitor the state of an already-stable system just
in case something goes wrong. Instead it is more efficient to only
sporadically check on the system to make sure things are behaving
well. The general topic of research on these types of problems is then
to determine precisely when control signals should be updated to
improve efficiency while still guaranteeing a desired quality of
service. Many different researchers have explored these ideas over the
past five decades under many different names including `event-based
sampling', `adaptive sampling', `event-triggered control', `Lebesgue
sampling', `send-on-delta concept', or `minimum attention
control'~\citep{KJA-BMB:99,KJA-BMB:02,KEA:99,JS-WPMHH-PB:05,RWB:97,MM:06,WPMHH-JHS-PPJB:08,PT:07,PGO-JRM-DMT:02,JM-WM:67,AS-CP-NM:14}. While
these ideas have been tossed around since the 1960's, it is only in
the last 10 years that the field has started maturing to soon stand
alone in the area of systems and
control~\citep{WPMHH-KHJ-PT:12,LH-CF-HO-AS-EF-JR-SIN:17}, with a
specific set of interesting challenges in the context of network
systems.

\subsection*{Why Event-Triggering?}

The idea of event-triggered sampling and control stems from an
implementation issue; how can we control continuous-time systems using
digital controllers? The current standard method is to simply have the
digital controller take actions periodically; and while ideas of
aperiodic sampling and control have been proposed long ago, some
modern textbooks seem to suggest that periodic sampling and control is
the only way to implement feedback control laws on digital
systems~\citep{KJA-BW:97,GF-JP-AE:10,WPMHH-KHJ-PT:12}.  However, real
systems in general are not able to acquire samples at an exact
operating frequency. Consequently, the stability of sampled-data
systems with aperiodic sampling has been a longstanding research area
within the systems and controls community.  As a result, there are
already many different standard methods and ways of thinking about
such problems and analyzing stability. For example, aperiodic sampling
can actually just be modeled as a specific time-delay system. The same
system might be modeled as a hybrid system with impulsive
dynamics. More specifically, a Linear Time Invariant (LTI) system with
aperiodic sampled-data might be reformulated as a Linear Time Varying
(LTV) or Linear Parameter Varying (LPV) system. Another option is to
derive the input/output relationships to study the effect of aperiodic
sampling on an output as is often done in robust control. In any case,
the main question of interest in a nutshell is: how quickly does the
system need to be sampled to guarantee stability? The recent
survey~\citep{LH-CF-HO-AS-EF-JR-SIN:17} presents various methods to
address this question.

\subsubsection*{Aperiodic Sampling as an Opportunity}

Whereas the above paradigm views aperiodic sampling as a type of
disturbance with respect to the ideal case of exactly periodic
sampling, new advances in event-triggering methods suggest treating
aperiodic sampling as an opportunity rather than an inconvenience or
disturbance. As mentioned above, these types of problems generally
arise as we try to control continuous-time systems with digital
controllers using an idealized controller that assumes exact state
information and continuous feedback is always possible. The natural
question to ask at this point is then exactly how fast does the
controller need to sample the system and feed back the control input
to ensure closed-loop stability? The answers to these questions then
often come in terms of robustness guarantees to the tune of ``as long
as the sampling/control frequency is greater than some threshold, then
the steady-state error is guaranteed to remain less than some
quantity.'' In other words, as long as the feedback loop for a given
system is `fast enough', the system behaves similarly to the ideal
system. In particular, this paradigm hinges largely on the intuitive
fact that, as long as the sampling rate is sufficiently fast, the
system behaves well. However, early studies into this question
revealed that this is not always the case, e.g., it is actually
possible to speed up the sampling rate and have a closed-loop feedback
control system degrade in
performance~\citep{SG:63,GB-RT:66,AL-JW:66,RT-GB:66,JM-WM:67,DC-LM:67}. However,
this was a non-issue in general due to the fact that while it is not
guaranteed that speeding up the sampling rate improved performance at
all times, it is true that continuing to increase the sampling speed
will eventually yield better performance at some point. While this
paradigm has been mostly sufficient for controlling many different
autonomous systems in the past, it seems quite limiting in many
application areas today. Instead, event-triggered ideas essentially
recast the question of `how fast should a control system respond' as
`exactly when should a control system respond?' to improve efficiency.

\subsection*{Why Now?}

When considering a dedicated sensor and actuator that are not
connected to any wireless network, it may be reasonable to ask the
sensor to take samples as fast as possible at all times and have the
actuator act accordingly. In this setting, it is probably not even
worth developing a more efficient or intelligent sensor as the
dedicated sensor periodically taking measurements is not affecting
anything else. However, this may not be practical if we only have a
remote sensor and the sensor data must be transmitted back to the
actuator over a wireless channel, especially if this wireless channel
must be shared among other devices.  In this case the action of
acquiring a sample is now literally a resource that must be managed
efficiently.  More generally, a strong motivation for the resurgence
of these topics is likely due to the increasing popularity of
networked cyber-physical systems across all disciplines. In
particular, the inherently tight couplings required between physical
processes (e.g., sampling, actuation, motion) and cyber processes
(e.g., communication, computation, storage) in networked systems
reveals the need for more efficient deployment of such systems by
treating things like wireless communication or computation as
resources rather than taking them for granted. This suggests the
application of event-triggered ideas not only to determine when
control signals should be updated, but to a wider array of
capabilities including data acquisition, actuation, and communication.
{It is in this sense that, hereafter, we employ the term
  event-triggered ``communication'' to refer to a communication event
  and the term event-triggered ``control'' to refer to a controller
  update event. When both appear in conjunction, we refer to the
  combination of event-triggered communication and control as
  event-triggered ``coordination''.}

In particular, we focus on the Internet of Things (IoT) and other
large-scale networks as a strong motivator for why we should think
about event-triggered coordination schemes rather than periodic ones.
IoT devices need to support a large variety of sensors and actuators
that interact with the physical world, in addition to standard cyber
capabilities such as processing, storage, or communication.  However,
as IoT devices aim to support services and applications that interact
with the physical world,  large numbers of these devices need to
be deployed and work reliably with minimal human
intervention~\citep{CP-CHL-SJ-MC:14}. This requirement places a lot of
crucial constraints on what we expect of our IoT
devices~\citep{PK-CP-GE-MP:16}. First of all, these devices in general
will be battery-operated and have small form factors, making energy
efficiency a critical design consideration. Second, these devices will
need to have a wide range of capabilities to integrate seamlessly
within a larger IoT network, translating to high computational
complexities. Third, the majority of communications within IoT
networks are wireless, meaning wireless congestion is another
important consideration. Consequently, the cyber operations
(processing, storage, and communication) can no longer be taken for
granted and must instead be viewed as a scarce, globally shared
resource.

Given these new-age design considerations, it is not surprising that
event-triggering has recently been gaining a lot of traction as a
promising paradigm for addressing the issues
above~\citep{CP-AZ-PC-DG:14,PK-CP-GE-MP:16}. Event-triggered methods
are useful here in that they address precisely \emph{when} different
actions (e.g., sensor sample, wireless communication) should occur to
efficiently maintain some desired property. The resurgence we see now
may be credited to the seminal works~\citep{KJA-BMB:99,KJA-BMB:02,KEA:99}, where
the advantages of event-triggered control over periodic
implementations were highlighted. Interestingly, \cite{AA-RB:04}
compared the differences between event-triggered and time-triggered
distributed control systems and concluded that one of the main
deficiencies of periodic control is in its lack of flexibility and
scalability. Given the current vision of the IoT being extremely
massive and interactive, it is clear we need methods to help enhance
flexibility and scalability for these systems of the future. These
ideas have thus been gaining more momentum which we also credit partly
to the rise of networked control systems in general.

\subsection*{Technical Challenges {Specific to Networked}
  Systems}

Early works on the subject assume a single decision-maker is
responsible for when different actions should be taken by a system;
however, we now need ways of implementing these ideas in fully
distributed settings to be applicable to the
IoT. \cite{WPMHH-KHJ-PT:12} provide a survey of event-triggered and
self-triggered control but focus on scenarios where the events are
being dictated by a single decision-maker or controller. Instead, the
focus of this article is on how to extend these ideas to distributed
settings and the new technical challenges that must be addressed in
doing so. In particular, we must emphasize the fact that applying the
ideas of event-triggered control to networked systems in general poses
many new challenges that do not exist in either area alone. For
instance, event-triggered coordination algorithms automatically
introduce asynchronism into a system which makes their analysis more
difficult. Furthermore, it often becomes difficult to find local
triggering rules that agents with distributed information can apply to
ensure some system level properties are satisfied; whereas in a
centralized setup it is generally easier to find a triggering rule
that can directly control some quantity of interest. For example, it
is easy to constrain a centralized decision maker to allocate at most
a certain number of actions per time period; however, it is more
difficult to distribute these decisions to both be efficient and still
be sure that the total number of actions per time period constraint is
respected.

\subsection*{Why Consensus?}

% While we may like to survey the field of applying
% event/self-triggered control ideas to networked systems in general,
% doing so would likely not be possible nor useful. As a result,

Given the wide variety of opportunistic state-triggered control ideas
in networked systems, we have made a conscious decision to focus
specifically on consensus problems, as a canonical example of
distributed algorithms in general. Nevertheless, our discussion
illustrates many of the challenges that arise beyond the specific
problem of consensus.  For example, it is already known that the
separation principle does not hold for event-triggered control systems
in general~\citep{CR-HS-LB-KHJ:11}.  Since the idea of event-triggered
coordination is to take various actions when only deemed necessary,
the specific task at hand is tightly coupled with when events should
be triggered. However, event-triggered algorithms are certainly not
unique solutions to any given problem either. Given a specific problem
instance, there are many different event-triggered algorithms that can
solve the problem.  By choosing a simple, but concrete set of
problems, this article discusses many different event-triggered
algorithms that have been recently proposed and what exactly the
seemingly subtle differences are. At the same time, the set of
consensus problems is still general enough that the methods/reasoning
behind the event-triggered algorithms we discuss throughout this
article are applicable to a number of different application areas
related to networked systems. The problem set is also rich enough to
capture the same technical difficulties that arise in many other
networked event-triggered scenarios such as how to deal with the
natural asynchronism introduced into the systems and how to guarantee
Zeno behaviors are excluded (e.g., algorithm certificates including
deadlocks being avoided).

% By focusing specifically on consensus problems, we present a number of
% different event-triggered coordination algorithms and discuss in
% detail their derivations, discuss the motivations behind different
% triggering rules, and highlight the technical challenges behind their
% design and analysis.

\subsection*{Organization}

We begin in Section~\ref{se:background} by formulating the basic
multi-agent consensus problem and provide a short background on the
first time event-triggered ideas were applied to it. We then close the
section by identifying five different categories of properties to help
classify different pairs of consensus problems and event-triggered
coordination solutions. In Section~\ref{se:main}, we provide details
behind the five different categories including the motivations behind
them, and provide numerous examples of different algorithms that fall
under the different classifications. Proofs of most results discussed
in the article are presented in the Appendix.  In
Section~\ref{se:applications}, we take a step back from consensus by
providing many different general networking areas that can both
directly and indirectly benefit from the ideas discussed in this
article.  {In Section~\ref{se:beyond}, we provide an outlook
  on the role of event-triggered coordination in networked systems
  beyond multi-agent consensus and discuss interesting future lines of
  research.}  Finally, we gather some concluding thoughts in
Section~\ref{se:conclusions}.

\subsection*{Preliminaries}

We introduce some notational conventions used throughout the
article. Let $\real$, $\realpositive$, $\realnonnegative$, and
$\integerpositive$ denote the set of real, positive real, nonnegative
real, and positive integer numbers, respectively. We denote by
$\ones{N}$ and $\zeros{N} \in \real^N$ the column vectors with entries
all equal to one and zero, respectively. The~$N$-dimensional identity
matrix is denoted by~$I_N$. Given two
matrices~$A \in \real^{m \times n}$ and $B \in \real^{p \times q}$, we
denote by $A \otimes B \in \real^{mp \times nq}$ as their Kronecker
product. We let $\Norm{\cdot}$ denote the Euclidean norm on
$\real^N$. We let
$\diag{\real^N} = \setdef{x \in \real^N}{x_1=\dots=x_N} \subset
\real^N$
be the agreement subspace in $\real^N$.  For a finite set $S$, we let
$|S|$ denote its cardinality. Given $x, y \in \real$, Young's
inequality states that, for any $\epsilon \in \realpositive$,
\begin{align}\label{eq:Young}
  xy \leq \frac{x^2}{2\epsilon} + \frac{ \epsilon y^2}{2}.
\end{align}
A weighted directed graph (or weighted digraph) $\commgraph = (V,E,W)$
is comprised of a set of vertices $V = \until{N}$, directed edges
$E \subset V \times V$ and weighted adjacency matrix
$W \in \realnonnegative^{N \times N}$. Given an edge $(i,j) \in E$, we
refer to $j$ as an out-neighbor of $i$ and $i$ as an in-neighbor of
$j$.  The weighted adjacency matrix $W \in \real^{N \times N}$
satisfies $w_{ij} > 0$ if $(i, j) \in E$ and $w_{ij} = 0$
otherwise. The sets of out- and in-neighbors of a given node $i$ are
$\Nouti$ and $\Nini$, respectively. The graph~$\commgraph$ is
\emph{undirected} if and only if $w_{ij} = w_{ji}$ for all
$i,j \in \vertices$. For convenience, we denote the set of neighbors
of a given node~$i$ in an undirected graph as simply~$\NN_i$.  A path
from vertex $i$ to $j$ is an ordered sequence of vertices such that
each intermediate pair of vertices is an edge.  A digraph $\commgraph$
is strongly connected if there exists a path between any two vertices.
The out- and in-degree matrices $\Dout$ and $\Din$ are diagonal
matrices where
\begin{align*}
  \douti = \sum_{j \in \Nouti} w_{ij} , \quad d_i^{\operatorname{in}} =
  \sum_{j \in \Nini} w_{ji} ,
\end{align*}
respectively. A digraph is weight-balanced if $\Dout = \Din$.  The
(weighted) Laplacian matrix is $L = \Dout - W$. Based on the structure
of $L$, at least one of its eigenvalues is zero and the rest of them
have nonnegative real parts.  If the digraph $\commgraph$ is strongly
connected, $0$ is a simple eigenvalue with associated eigenvector
$\ones{N}$.  The digraph~$\commgraph$ is weight-balanced if and only
if $\ones{N}^T L = \zeros{N}$, which is also equivalent to
$\Sym{L}=\frac{1}{2}(L + L^T)$ being positive semidefinite. For a
strongly connected and weight-balanced digraph, zero is a simple
eigenvalue of $\Sym{L} $.  In this case, we order its eigenvalues as
$\lambda_1 = 0<\lambda_2\leq \dots \leq \lambda_N$, and note the
inequality
\begin{align}\label{eq:LapBound}
  % \lambda_N(\Sym{L}) \| x -\frac{1}{n} (\ones{N}^T x) \ones{N} \|^2
  % \ge
  x^T L x \ge \lambda_2(\Sym{L}) \| x -\frac{1}{N} (\ones{N}^T x)
  \ones{N} \|^2 ,
\end{align}
for all $x \in \real^N$.  The following property will also be of use
later,
\begin{align}\label{eq:lap2-bound}
  \lambda_2(\Sym{L}) x^T L x \leq x^T \Sym{L}^2 x \leq
  \lambda_N(\Sym{L}) x^T L x .
\end{align}
This can be seen by noting that $\Sym{L}$ is diagonalizable and
rewriting $\Sym{L} = S^{-1} D S$, where $D$ is a diagonal matrix
containing the eigenvalues of~$\Sym{L}$.
% It can then be shown that
% \begin{align*}
%   \lambda_2(L) x^T S^{-1} D S x \leq x^T L^2 x \leq \lambda_N(L) x^t
%   S^{-1} D S x .
% \end{align*}

%{\color{green!60!black}
%  \cite{AA-FA-DL-GS-DVD-MB-KHJ:15} - leader follower, switching,
%  nonlinear synchronization - not directly applicable to average
%  consensus
%
%  \cite{GG-LD-QH:14} - linear systems $Ax_i + Bu_i$ - not directly
%  applicable to average consensus }

\section{What is Event-Triggered Consensus?}\label{se:background}

In this section we formally state the problem of event-triggered
consensus, which results from the application of event-triggered
control to the multi-agent consensus problem. We first describe the
basic approach to event-triggered control design and then
particularize our discussion to event-triggered consensus. 

\subsection{A primer on event-triggered control}\label{se:primer}

{We start by informally describing the event-triggered design approach
to stabilization along the lines proposed in~\citep{PT:07}}.  Given a
system on $\real^n$ of the form
\begin{align*}
  \dot x = F(x,u)
\end{align*}
with an unforced equilibrium at~$x^*$ (i.e., $F(x^*,0) = 0$), the
starting point is the availability of
\begin{enumerate}
\item a continuous-time controller $\map{k}{\real^n}{\real^m}$, along
  with
\item a certificate of its correctness in the form of a Lyapunov
  function $\map{V}{\real^n}{\real}$.
\end{enumerate}
In other words, the closed-loop system $\dot x = F(x,k(x))$ makes
$x^*$ asymptotically stable, and this fact can be guaranteed through
$V$ as Lyapunov function. The idea of event-triggered control is the
following: rather than continuously updating the input $u$ as $k(x)$,
use instead a sampled version $\hat{x}$ of the state to do it as
$k(\hat{x})$. This sample gets updated in an opportunistic fashion in
a way that still ensures that $V$ acts as a certificate of the
resulting sampled implementation. If done properly, this has the
advantage of not requiring continuous updates of the input while still
guaranteeing the original stabilization of the equilibrium point. The
question is then how to determine when the sampled state needs
updating.  Formally, the closed-loop dynamics looks like
\begin{align}\label{eq:dynamicsgeneral}
  \dot x = F(x,k(\hat{x}))  
\end{align}
and hence one has $\dot V = \nabla V(x) \cdot F(x,k(\hat{x}))$.
{More specifically, letting $\{ t_\ell \}_{\ell \in
    \integernonnegative}$ denote the sequence of event times at which
  the control input is updated, the control input is given by
\begin{align}\label{eq:control}
  u(t) &= k(\hat{x}(t)),
\end{align}
where the sampled state is given by
\begin{align}\label{eq:estimate}
\hat{x}(t) &= x(t_\ell) \quad \text{for} \quad t \in [t_\ell,t_{\ell+1}),
\end{align}
for some sequence of time~$\{t_\ell\}_{\ell \in \integernonnegative}$. 
In other words, the control signal~$u(t)$ is only updated at the
discrete times~$t_\ell$ and the input is held constant in between
events. The goal is to determine a specific event-condition such that
the closed-loop system still converges to the desired state.

Under mild conditions on~$F, k,$ and~$V$ (formally, $F$ uniformly -in
$x$- Lispchitz in its second argument, $k$ Lipschitz, and $\nabla V$
bounded), some manipulations of~$\dot{V}$ leads to expressions of the
form
\begin{align}\label{eq:auxxx}
  \dot V \le \nabla V(x) \cdot F(x,k(x)) + G(x) \|e\| ,
\end{align}
for some function $G$ taking nonnegative values, and where $e =
\hat{x} - x$ is the error between the sampled and the actual state.}
The first term of the derivative in~\eqref{eq:auxxx} is negative,
while the second vanishes when the sampled state coincides with the
actual one, i.e., $e=0$. Therefore, to ensure that $\dot V<0$, one can
simply define a trigger that prescribes that the sampled state should
be updated whenever the magnitudes of the first and second term are
equal.

{This is encoded through what is called a \emph{triggering
    function} or \emph{event-trigger}~$f(\cdot)$, which evaluates whether a given
  state~$x$ and error~$e$ combination should trigger an event or
  not. With a slight abuse of notation, we define this condition as
\begin{align}\label{eq:generaltrigger}
  f(e,w) \triangleq g(e) - h(w) = 0,
\end{align}
where~$g : \real^n \rightarrow \realnonnegative$ is a nonnegative
function of the error with~$g(0) = 0$ and~$h \in \realnonnegative$ is
a threshold function that may depend on variables like the state~$x$,
the sampled state~$\hat{x}$, or time~$t$, and even additional
variables or parameters. For now, we lump them all together in the
variable $w$, and as we make progress in our exposition, we detail
what $w$ is in each case.  The point of this triggering
function~\eqref{eq:generaltrigger} is that it guarantees some function
of the error~$g(e)$ is always smaller than some
threshold~$h(\cdot)$. This happens because when the
condition~\eqref{eq:generaltrigger} is satisfied, an event is
triggered, which resets the error~$e = 0$ and thus~$g(e) = 0$ is also
reset.  Specifically, given a triggering function, event times are
implicitly defined  by
\begin{align}\label{eq:eventtimes}
  t_{\ell+1} = \min \{ t' \geq t_\ell \hspace*{1mm} | \hspace*{1mm}
  f(e(t'),w(t')) = 0 \} .
\end{align}

We are then interested in designing these functions~$g$ and~$h$ in
such a way that the closed-loop dynamics~\eqref{eq:dynamicsgeneral}
with control inputs~\eqref{eq:control} driven by~\eqref{eq:eventtimes}
ensures~$x \rightarrow x^*$. Using~\eqref{eq:auxxx}, it is easy to see
that the triggering function with $w=x$ defined by
\begin{align*}
  g(e) &= \TwoNorm{e},
  \\
  h(x) &= \frac{|\nabla V(x) \cdot F(x,k(x))|}{|G(x)|},
\end{align*}
ensures that
\begin{align*}
  \dot V \le \nabla V(x) \cdot F(x,k(x)) + G(x) \|e\| < 0
\end{align*}
at all times. This fact can ultimately be used to show that $x$
asymptotically approaches~$x^*$ as long as there are no deadlocks in
the execution. We discuss this point in detail next.
}

\subsubsection{Deadlocks, Zeno behavior, and \changes{Minimum
    Inter-Event Time}}\label{se:zeno}

With the trigger design in place, one can analyze the behavior of
the resulting implementation, such as guaranteeing liveness and the
absence of deadlocks. We formalize this in the following definition.

\begin{definition}[Zeno behavior]\label{def:zeno}
  Given the closed-loop dynamics~\eqref{eq:dynamicsgeneral} with
  control inputs~\eqref{eq:control} driven by~\eqref{eq:eventtimes}
  % according to some triggering condition~$f(\cdot)$, we say 
  a solution with initial condition~$x(0) = x_0$ exhibits \emph{Zeno
    behavior} if there exists $T > 0$ such that $t_{\ell} \leq T$ for
  all $\ell \in \integernonnegative$.
\end{definition}

In other words, if the event-triggered controller defined by the
triggering function~\eqref{eq:generaltrigger} demands that an infinite
number of events (e.g., controller updates) occur in a finite time
period, the solution exhibits Zeno behavior. Note that it is possible
that depending on the initial condition~$x(0) \in \real^n$, different
solutions may or may not exhibit Zeno behavior. Only in the case when
it is guaranteed that Zeno behavior does not occur along \emph{any}
trajectory, we say that the system (as a whole) does not exhibit Zeno
behavior.

\changes{Being able to rule out Zeno behavior is extremely important
  in validating the correctness of a given event-triggered controller.
  In general, the event-triggered algorithms we discuss are comprised
  of some kind of control law and triggering rule, with the latter
  driving what information is being used by the control law in real
  time. The existence of Zeno behavior means there exists an
  accumulation time $T > 0$ by which an infinite number of events will
  be triggered. This is problematic for any physical implementation on
  a real-time system, as it is asking the controller to be updated
  with new information an infinite number of times in a finite time
  period.

  Another point worth highlighting is the difference between ruling
  out Zeno behavior versus ensuring a uni- form minimum time between
  any two consecutive events. In fact, the guarantee on lack of Zeno
  behavior is weaker than ensuring that there exists a quantity
  $\tau^\text{min}$ that uniformly lower bounds the time in between
  consecutive events, i.e.,
  \begin{align*}
    t_{\ell+1} - t_\ell \geq \tau^\text{min} > 0
  \end{align*}
  for all $\ell \in \integernonnegative$, which is a more pragmatic
  property when considering physical hardware. We refer
  to~$\tau^\text{min}$ as the minimum inter-event time
  (MIET)~\citep{DPB-WPMHH:14}.  Since dedicated hardware can only
  operate at some maximum frequency (e.g., a physical device can only
  broadcast a message or evaluate a function a finite number of times
  in any finite period of time), ensuring the existence of a positive
  MIET is more appropriate for physical implementation that simply
  ruling out Zeno behavior.
  % This means that simply ensuring a system does not exhibit
  % Zeno behavior may not be enough to guarantee the algorithm can be
  % implemented on a physical system, if the physical hardware cannot
  % match the speed of actions required by the algorithm. Thus to
  % guarantee not only theoretical correctness of an event-triggered
  % algorithm, it is also imperative to ensure the existence of a
  % positive MIET in order to guarantee that the solution can
  % physically
  % be implemented on a real system.

  Next, we provide examples describing the seemingly subtle
  differences between these concepts; and more importantly, their
  implications on correctness and implementation. Consider the dynamic
  system~\eqref{eq:dynamicsgeneral} for which a triggering
  function~$f$ has been defined as in~\eqref{eq:generaltrigger} that
  leads to three different sequences of event times~$\{t_\ell\}_{\ell \in
    \integernonnegative}$ described by~\eqref{eq:eventtimes}:
\begin{description}
\item[(Zeno behavior):] Consider
  \begin{align*}
    t_{\ell+1}-t_\ell = \frac{1}{(\ell+1)^2},
  \end{align*}
  for~$\ell \in \integernonnegative$. Given~$t_0 = 0$, this defines
  the sequence of times as
  \begin{align*}
    t_\ell = \sum_{n=1}^\ell \frac{1}{n^2}.
  \end{align*}
  As the number of events~$\ell \rightarrow \infty$, we have that
  $t_\ell \leq \frac{\pi^2}{6}$ for all~$\ell \in
  \integernonnegative$. This means that even if there existed a
  physical device that can perform actions this quickly, the
  theoretical analysis of the closed-loop dynamic
  system~\eqref{eq:dynamicsgeneral} is not valid beyond~$T =
  \frac{\pi^2}{6}$.
\item[(Non-Zeno behavior without a positive MIET):] Consider
  \begin{align*}
    t_{\ell+1}-t_\ell = \frac{1}{\ell+1},
  \end{align*}
  for~$\ell \in \integernonnegative$. Given~$t_0 = 0$, this defines
  the sequence of times as
  \begin{align*}
    t_\ell = \sum_{n=1}^\ell \frac{1}{n}.
  \end{align*}
  In this case, as~$\ell \rightarrow \infty$ we also have that~$t_\ell
  \rightarrow \infty$, which means Zeno behavior can be
  excluded. However, since the inter-event times~$t_{\ell+1} -
  t_{\ell}$ go to 0 as $\ell \rightarrow \infty$, there does not exist
  a positive MIET~$\tau^\text{min}$. This means that although the
  theoretical analysis might guarantee stability of the closed-loop
  dynamic system~\eqref{eq:dynamicsgeneral}, it would require hardware
  that can perform actions infinitely fast.
\item[(Positive MIET):] Consider
  \begin{align*}
    t_{\ell+1} - t_\ell = c + \frac{1}{\ell+1},
  \end{align*}
  for some~$c > 0$ and all~$\ell \in \integernonnegative$. Given~$t_0
  = 0$, this defines the sequence of times as
  \begin{align*}
    t_\ell = \sum_{n=1}^\ell cn + \frac{1}{n}.
  \end{align*}
  We can now guarantee not only the absence of Zeno behavior, but that
  there exists a positive MIET such that all inter-event times are
  lower-bounded $t_{\ell+1} - t_{\ell} \geq \tau^\text{min} = c >
  0$. This not only guarantees stability of the closed-loop dynamic
  system~\eqref{eq:dynamicsgeneral}, but also that the solution can
  actually be implemented using a device that can take actions at a
  frequency faster than~$\frac{1}{\tau^\text{min}}$.
\end{description}

Based on the above discussion, it is important to realize that a
complete, fully implementable event-triggered control solution to a
problem should also include the existence of a positive
MIET.
% Unfortunately, there are currently many
% published works that suggest event-triggered solutions to various
% problems that do not satisfy all of these properties, which leaves
% questions about how practical they are.

}

\subsection{Multi-agent average consensus}

Here, taking as reference our discussion above, we proceed to describe
the multi-agent average consensus problem, identifying as we go the
key elements (continuous-time controller and certificate) necessary to
tackle the design of event-triggered coordination mechanisms.  We
start with a simple, yet illustrative, scenario to introduce the main
ideas. Towards the end of the section, we discuss various directions
along which the problem and its treatment gains in complexity and
realism.  

We let $\commgraph$ denote the connected, undirected graph
that describes the communication topology in a network of~$N$
agents. In other words, agent~$j$ can communicate with agent~$i$ if
$j$ is a neighbor of~$i$ in~$\commgraph$.  We denote by
$x_i \in \real$ the state of agent $i \in \until{N}$ and consider
single-integrator dynamics
\begin{align}\label{eq:dynamics}
  \dot{x}_i(t) = u_i(t).
\end{align}
It is well-known that the distributed controller
\begin{align}\label{eq:ideal}
  u_i^*(t) = -\sum_{j \in \NN_i} \left( x_i(t) - x_j(t) \right)
\end{align}
drives the states of all agents to the average of the initial
conditions~\citep{ROS-RMM:03c,ROS-JAF-RMM:07}. This is formalized in
Theorem~\ref{th:main}.

\begin{theorem}\longthmtitle{Continuous
    controller~\citep{ROS-RMM:03c}}\label{th:main}
  Given {a connected, undirected graph~$\commgraph$ and} the
  dynamics~\eqref{eq:dynamics}, if all agents implement the control
  law~\eqref{eq:ideal}, then the system asymptotically achieves
  multi-agent average consensus; i.e.,
  \begin{align}\label{eq:averageconsensus}
    \lim_{t \rightarrow \infty} x_i(t) = \frac{1}{N} \sum_{j=1}^N
    x_j(0)
  \end{align}
  for all $i \in \until{N}$.
\end{theorem}
%\begin{pf}
%In the appendix. \hfill \qed
%\end{pf}

Implementing~\eqref{eq:ideal} in a digital setting is not possible
since it requires all agents to have continuous access to the state of
their neighbors and the control inputs $u_i(t)$ must also be updated
continuously. This is especially troublesome in the context of
wireless networked systems since this means agents must communicate with
each other continuously as well. Instead, researchers have been
interested in applying event-triggered strategies to relax these requirements.

\subsection{Centralized event-triggered control}

Consider the dynamics~\eqref{eq:dynamics} and the ideal control
law~\eqref{eq:ideal}.  Letting $x = (x_1, \dots, x_N)^T$ and
$u = (u_1, \dots, u_N)^T$, the closed-loop dynamics of the ideal
system is given by
\begin{align}\label{eq:idealfull}
  \dot{x}(t) = -Lx(t) ,
\end{align}
where~$L$ is the Laplacian of~$\commgraph$.  As stated before,
implementing this requires all agents to continuously update their
control signals which is not realistic for digital
controllers. Instead, let us consider a digital implementation of this
ideal controller
\begin{align}\label{eq:centralcontrol}
  u(t) = -Lx(t_\ell), \quad t \in [t_\ell, t_{\ell+1} ),
\end{align}
where the event times $\{ t_\ell \}_{\ell \in \integernonnegative}$
are to be determined such that the system still converges to the
desired state. It is worth mentioning here that the control
law~\eqref{eq:centralcontrol} is chosen such that the average of all
agent states is an invariant quantity regardless of how the event
times~$\{t_\ell\}_{\ell \in \integernonnegative}$ are chosen, thus
preserving the average of the initial conditions throughout the
evolution of the system.  More specifically, utilizing this
controller,
\begin{align}\label{eq:average-preservation}
  \frac{d}{dt} (\ones{N}^T x(t)) = \ones{N}^T \dot x(t) = \ones{N}^T L
  x(t_\ell) = 0 ,
\end{align}
where we have used the fact that $L$ is symmetric and
$L \ones{N} = 0$.

Let $e(t) = x(t_\ell) - x(t)$ for $t \in [t_\ell, t_{\ell+1})$ be the
state measurement error.  For simplicity, we denote by
$\hat{x}(t) = x(t_\ell)$ for $t \in [t_\ell, t_{\ell+1})$ as the state
that was used in the last update of the control signal.  The
closed-loop dynamics of the controller~\eqref{eq:centralcontrol} is
then given by
\begin{align}\label{eq:centraldynamics}
  \dot{x}(t) = -L \hat{x}(t) = -L( x(t) + e(t) ).
\end{align}

%\margin{should we fully formalize these problems? I think no...}
%{\color{red}
%\begin{problem}[Centralized event-triggered control]
%Given the closed-loop dynamics
%\begin{align*}
%\dot{x}(t) &= - L \hat{x}(t) = - L (x(t) + e(t)) , \\
%\dot{e}(t) &= - \dot{x}(t) = L (x(t) + e(t)), \\
%e(t_\ell)^+ &= 0, \quad \ell \in \integernonnegative, \\
%t_{\ell+1} &= \min\{ t' \geq t_\ell | f(x(t'),e(t')) = 0 \},
%\end{align*}
%find an event-trigger~$f(x,e)$ such that the sequence of times~$\{t_\ell\}_{\ell \in \integernonnegative}$ ensures multi-agent average consensus is achieved.
%\end{problem}
%}

{We are then interested in designing a triggering condition of the
form~\eqref{eq:generaltrigger} in such a way that the closed-loop
dynamics~\eqref{eq:centraldynamics} driven by~\eqref{eq:eventtimes}
ensures multi-agent average consensus is achieved. The problem can now
be formalized as follows.}

\begin{problem}\longthmtitle{Centralized event-triggered {consensus}}\label{pr:central}
  Given the closed-loop dynamics~\eqref{eq:centraldynamics}, find an
  event-trigger~$f(\cdot)$ such that the sequence of times~$\{ t_\ell
  \}_{\ell \in \integernonnegative}$ ensures multi-agent average
  consensus~\eqref{eq:averageconsensus} is achieved.
\end{problem}

Some of the first works to consider this problem were
\cite{DVD-EF:09,DVD-KHJ:09,EK-XW-NH:10}. Following~\citep{DVD-EF-KHJ:12},
to solve this problem we consider the Lyapunov function
\begin{align*}
  V(x) = \frac{1}{2}x^T L x.
\end{align*}
Given the closed-loop dynamics~\eqref{eq:centraldynamics}, we have
\begin{align}\label{eq:lyapderivative1}
  \dot{V} = x^T L \dot{x} = -x^T L L (x + e) =
  -\underbrace{\TwoNorm{Lx}^2}_\text{"good"} - \underbrace{x^T L L
    e}_\text{"bad"} .
\end{align}
The main idea of (Lyapunov-based) event-triggered control is then to
determine when the controller should be updated (i.e., when the error
$e$ should be reset to 0) by balancing the ``good'' term against the
potentially ``bad'' term. More specifically, we are interested in
finding conditions on the error~$e$ such that $\dot{V} < 0$ at all
times. Using norms, we can bound
\begin{align*}
  \dot{V} \leq -\TwoNorm{Lx}^2 + \TwoNorm{Lx}\TwoNorm{L}\TwoNorm{e}.
\end{align*}
Then, if we can somehow enforce the error~$e$ to satisfy
\begin{align*}
  \TwoNorm{e} \leq \sigma \frac{ \TwoNorm{Lx} }{\TwoNorm{L}} ,
\end{align*}
with $\sigma \in (0,1)$ for all times, we have
\begin{align*}
  \dot{V} \leq (\sigma - 1) \TwoNorm{Lx}^2,
\end{align*}
which is strictly negative for all $Lx \neq 0$. It is then easy to see
that the following centralized event-trigger using
\begin{align*}
  g(e) &= \TwoNorm{e} , 
  \\
  h(x) &= \sigma \frac{\TwoNorm{ Lx }}{\TwoNorm{L}} ,
\end{align*}
ensures this is satisfied at all times. Note that in this case we
have a state-dependent threshold~$h(x)$, but other types of thresholds
will be discussed later.

\begin{theorem}\longthmtitle{Centralized event-triggered
    control~\citep{DVD-EF-KHJ:12}}\label{th:central} 
  Given {a connected, undirected graph~$\commgraph$ and} the
  closed-loop dynamics~\eqref{eq:centraldynamics}, if the event times
  are determined as the times when
  \begin{align}\label{eq:centraltrigger}
    f(e,x) &\triangleq \TwoNorm{e} - \sigma \frac{
      \TwoNorm{Lx}}{\TwoNorm{L}} = 0,
  \end{align}
  then the system achieves multi-agent average consensus.
\end{theorem}

In other words, given a control update at time~$t_\ell$, the next time
$t_{\ell+1}$ the controller is updated is given by~\eqref{eq:eventtimes}, 
\begin{align*}
  t_{\ell+1} &= \min \setdef{t' \geq t_\ell}{\TwoNorm{e(t')} = \sigma
    \frac{\TwoNorm{Lx(t')}}{\TwoNorm{L}}} .%\\
     % t_{\ell+1} &= \min \setdef{t \geq t_\ell }{f(x(t),e(t)) = 0} ,
\end{align*}

%The algorithm is formalized in Table~\ref{tab:algorithm2}.
%
%\begin{table}[htb]
%  \centering
%  \framebox[.9\linewidth]{\parbox{.85\linewidth}{%
%      \parbox{\linewidth}{At times $t \in [t_\ell,t_{\ell+1})$, system (continuously) performs:}
%      \vspace*{-1.5ex}
%      \begin{algorithmic}[1]
%      \STATE set $\hat{x}(t) = x(t_\ell)$
%      \STATE set $e(t) = \hat{x}(t) - x(t)$
%        \IF{$\TwoNorm{e(t)} = \sigma \frac{\TwoNorm{L x(t)}}{\TwoNorm{L}}$}  
%        \STATE set $t_{\ell+1} = t$
%        \STATE set $\hat{x}(t) = x(t_{\ell+1})$
%        \STATE set $\ell = \ell + 1$
%        \ENDIF
%        \STATE set $u(t) = - L \hat{x}(t)$
%      \end{algorithmic}}}
%  \caption{Centralized event-triggered control.}\label{tab:algorithm2}
%\end{table}

The proof of convergence to the desired state then follows almost
directly from the proof of Theorem~\ref{th:main} and the fact that the
sum of all states is still an invariant
quantity. Furthermore,~\cite{DVD-EF-KHJ:12} are able to rule out the
existence of Zeno behavior, cf. Definition~\ref{def:zeno}, \changes{by showing
there exists a positive MIET
\begin{align*}
  \tau^\text{min} = \frac{\sigma}{\TwoNorm{L}(1+\sigma)}
\end{align*}
uniformly bounding the inter-event times, i.e.,
\begin{align}\label{eq:lower-bound}
  t_{\ell+1} - t_\ell \geq \tau^\text{min} > 0 \quad \text{for all} \quad \ell
  \in \integernonnegative. 
\end{align}
}
\changes{As discussed in Section~\ref{se:zeno}, the existence of the positive MIET guarantees that the design is implementable over physical platforms.
}

The centralized event-triggered controller~\eqref{eq:centralcontrol}
with triggering law~\eqref{eq:centraltrigger} relaxes the requirement
that agents need to continuously update their control signals;
however, it still requires the controller to have perfect state
information at all times to be able to evaluate the triggering
condition~$f(\cdot)$. Next, we provide a distributed solution instead of
a centralized one.

\subsection{Decentralized event-triggered control}\label{se:canonical}

In the previous section we presented a centralized event-triggered
control law to solve the multi-agent average consensus
problem. Unfortunately, implementing this requires a centralized
decision maker and requires all agents in the network to update their
control signals simultaneously. Given the nature of and motivation
behind consensus problems, this is the first requirement we want to
get rid of. Here we present in detail the first real problem of
interest concerning this article.

Following~\citep{DVD-EF-KHJ:12}, consider a distributed digital
implementation of the ideal controller~\eqref{eq:ideal}. In this case
we assume each agent~$i$ has its own sequence of event times
$\{ t_\ell^i \}_{\ell \in \integernonnegative}$.  At any given
time~$t$, let
\begin{align}
  \hat{x}_i(t) = x_i(t_\ell^i) \text{ for } t \in [t_\ell^i, t_{\ell+1}^i )
\end{align} 
be the state of agent~$i$ at its last update time. The distributed
event-triggered controller is then given by
\begin{align}\label{eq:decentralcontrol}
  u_i(t) = - \sum_{j \in \NN_i} (\hat{x}_i(t) - \hat{x}_j(t) ) .
\end{align}
It is important to note here that the latest updated state
$\hat{x}_j(t)$ of agent $j \in \NN_i$ appears in the control signal
for agent~$i$. This means that when an event is triggered by a
neighboring agent~$j$, agent~$i$ also updates its control signal
accordingly.  As in the centralized case, let
$e_i(t) = x_i(t_\ell^i) - x_i(t)$ be the state measurement error for
agent~$i$. Then, letting $\hat{x} = (\hat{x}_1, \dots, \hat{x}_N)^T$
and $e = (e_1, \dots, e_N)^T$, the closed-loop dynamics of the
controller~\eqref{eq:decentralcontrol} is given by
\begin{align}\label{eq:decentraldynamics} 
  \dot{x}(t) = - L \hat{x}(t) = - L ( x(t) + e(t) ).
\end{align}
{Parallel to the general case in~\eqref{eq:eventtimes}, an
  event-trigger~$f_i(\cdot)$ for agent~$i$ is a function that
  determines its sequence of event times $\{ t_\ell^i \}_{\ell \in
    \integernonnegative}$ via
\begin{align}\label{eq:eventtimes-decentralized}
  t_{\ell+1}^i = \min \{ t' \geq t_\ell^i \hspace*{1mm} | \hspace*{1mm}
  f_i(e(t'),w(t')) = 0 \} .
\end{align}
 The problem we seek to solve can now be formalized as follows.}

{
  \begin{problem}\longthmtitle{Decentralized event-triggered
      consensus}\label{pr:event-control}
    Given a connected, undirected graph~$\commgraph$ and the
    closed-loop dynamics~\eqref{eq:decentraldynamics}, find an
    event-trigger~$f_i(\cdot)$ for each agent~$i$ that is locally
    computable and such that the sequences of times $\{ t_\ell^i
    \}_{\ell \in \integernonnegative}$ ensures multi-agent average
    consensus~\eqref{eq:averageconsensus} is achieved.
\end{problem}

By \emph{locally computable} function $f_i$, we mean that its value
only depends on variables that correspond to agent~$i$ and its
neighbors. Formally, this means that one can write
\begin{align*}
  f_i(e,w) = f_i(e_i,w_i) \triangleq g_i(e_i) - h_i(w_i) ,
\end{align*}
where $w_i$ represents information that is locally available to
agent~$i$.  Unlike in Problem~\ref{pr:central}, where we seek a single
event-triggering function~$f(\cdot)$ that depends on the global
state~$x$ to determine a global schedule, here we are interested in
having each agent~$i$ determine in a distributed way when its local
error~$e_i$ should be reset to 0.}

Following~\cite{DVD-EF-KHJ:12}, to solve this problem we again
consider the Lyapunov function
\begin{align*}
  V(x) = \frac{1}{2}x^T L x.
\end{align*}
Given the closed-loop dynamics~\eqref{eq:decentraldynamics}, we have
\begin{align}\label{eq:lyapderivative2}
  \dot{V} = -\TwoNorm{Lx}^2 - x^T L L e .
\end{align}
As before, we are interested in finding conditions on the error~$e$
such that~$\dot{V} < 0$ at all times; however, we must now do this in
a distributed way. For simplicity, let
$Lx \triangleq z = (z_1, \dots, z_N)^T$. Then, expanding out $\dot{V}$
yields
\begin{align*}
  \dot{V} &= - \left[ \sum_{i=1}^N z_i^2 - \sum_{j \in \NN_i} z_i (e_i - e_j) \right] \\
  &= - \left[ \sum_{i=1}^N z_i^2 - | \NN_i | z_i e_i + \sum_{j \in \NN_i} z_i
  e_j \right]. 
\end{align*}
Using Young's inequality~\eqref{eq:Young} and the fact
that~$\commgraph$ is symmetric, we can bound this by
\begin{align}\label{eq:finalbound}
  \dot{V} \leq - \left[ \sum_{i=1}^N (1 - a|\NN_i|)z_i^2 + \frac{1}{a}
  |\NN_i| e_i^2 \right]
\end{align}
for all $a > 0$. Letting $a \in (0, 1/|\NN_i|)$ for all $i$, if we can
enforce the error of all agents to satisfy
\begin{align*}
  e_i^2 \leq \frac{ \sigma_i a (1 - a|\NN_i|)}{|\NN_i|} z_i^2
\end{align*}
with $\sigma_i \in (0,1)$ for all times, we have
\begin{align}\label{eq:finalderivative}
  \dot{V} \leq \sum_{i=1}^N (\sigma_i - 1)(1-a|\NN_i|)z_i^2,
\end{align}
which is strictly negative for all $Lx \neq 0$. 
{In order to compute~$z_i$, agent~$i$ needs access to its
own state and its neighbors states,
\begin{align*}
  w_i = x_{\NN_i} \triangleq (x_i, \{x_j \}_{j \in \NN_i}).
\end{align*}
}
The following decentralized event-trigger then ensures that~$\dot{V}$ is strictly negative
until consensus has been achieved.

\begin{theorem}\longthmtitle{Decentralized event-triggered
    control~\citep{DVD-EF-KHJ:12}}\label{th:event-control}
  Given {a connected, undirected graph~$\commgraph$ and} the closed-loop
  dynamics~\eqref{eq:decentraldynamics}, if the event times of each
  agent~$i$ are determined as the times when
  % \begin{align}\label{eq:decentraltrigger}
  %   e_i^2 \triangleq \frac{ \sigma_i a (1 - a|\NN_i|)}{|\NN_i|}
  %   z_i^2,
  % \end{align}
  \begin{align}\label{eq:decentraltrigger}
    f_i(e_i,x_{\NN_i} ) \triangleq e_i^2 - \frac{
      \sigma_i a (1 - a|\NN_i|)}{|\NN_i|} z_i^2 = 0,
  \end{align}
  with $0 < a < 1/|\NN_i|$ for all~$i \in \until{N}$, then all
  non-Zeno trajectories of the system asymptotically achieve
  multi-agent average consensus.
\end{theorem}

Note that the trigger~\eqref{eq:decentraltrigger} can be evaluated by
agent~$i$ using only information about its own and neighbors' states.
{However, it should also be noted that implementing this algorithm
requires each agent~$i$ to have exact, continuous state information
about its neighbors~$\{x_j(t)\}_{j \in \NN_i}$. We address this in
Section~\ref{se:triggertype} below.}

%The algorithm is formalized in Table~\ref{tab:algorithm3}. 	
%
%\begin{table}[htb]
%  \centering
%  \framebox[.9\linewidth]{\parbox{.85\linewidth}{%
%      \parbox{\linewidth}{At times $t \in [t_\ell^i,t_{\ell+1}^i)$, agent~$i$ (continuously) performs:}
%      \vspace*{-1.5ex}
%      \begin{algorithmic}[1]
%      \STATE set $z_i(t) = \sum_{j \in \NN_i} (x_i(t) - x_j(t))$
%      %\STATE set $\hat{x}_i(t) = x_i(t_\ell^i)$
%      \STATE set $e_i(t) = \hat{x}_i(t) - x_i(t)$
%        \IF{$e_i(t)^2 = \frac{\sigma_i a (1 - a | \NN_i |)}{| \NN_i |} z_i(t)^2$}  
%        \STATE set $t_{\ell+1}^i = t$
%        \STATE broadcast $\hat{x}_i(t) = x_i(t_{\ell+1}^i)$ to neighbors $j \in \NN_i$
%        \STATE set $\ell = \ell + 1$
%        \ENDIF
%      \STATE set $u_i(t) = - \sum_{j \in \NN_i} (\hat{x}_i(t) - \hat{x}_j(t))$
%      \end{algorithmic}}}
%  \caption{Decentralized event-triggered control.}\label{tab:algorithm3}
%\end{table}

The proof of convergence to the desired state then directly follows
from the proof of Theorem~\ref{th:main} and the fact that the sum of
all states is still an invariant quantity. However, it is 
important to note that this argument is only valid along non-Zeno
trajectories, as discussed in Section~\ref{se:zeno}. Recall that the
result of Theorem~\ref{th:central} claimed all trajectories of the
system achieves multi-agent average consensus, but this was only
possible since it was already established in~\eqref{eq:lower-bound}
that Zeno behavior is impossible using the
trigger~\eqref{eq:centraltrigger} proposed in
Theorem~\ref{th:central} \changes{due to the existence of the positive MIET~$\tau^\text{min}$.}

Instead, in the derivation of the result of
Theorem~\ref{th:event-control},~\cite{DVD-EF-KHJ:12} only show that at
all times there exists one agent~$i$ for which the inter-event times
are strictly positive. Unfortunately, this is not enough to rule out
Zeno behavior, which is quite problematic, both from a pragmatic and
theoretical viewpoint, as the trajectories of the system are no longer
well-defined beyond the accumulation point in time. Consequently, the
main convergence result can only be concluded for trajectories that do
not exhibit Zeno behavior. Since Zeno behavior has in fact not yet
been ruled out for all trajectories using the
trigger~\eqref{eq:decentraltrigger}, the milder result of
Theorem~\ref{th:event-control} is all one can state.

{The intuitive reason for this is actually quite simple but it leads to
troubling implications: The main idea behind event-triggered control
is to only take certain actions when necessary. Since we are
interested in decentralized control protocols to achieve consensus for
a large system, it is easy to imagine some rare cases where some
agent~$i^*$ is already in agreement with its neighbors~$j \in
\NN_{i^*}$, but the rest of the system has not yet finished
evolving. In this case, once agent~$i$ reaches local consensus with
its neighbors, it wants to remain there. Unfortunately, this means
that the instant any of its neighbors begins to change its state
(because the rest of the network has not yet stabilized), the trigger
prescribes that agent~$i^*$ acts in response.

More specifically, looking at the trigger~\eqref{eq:decentraltrigger}
reveals that when~$z_{i^*} = (Lx)_{i^*} = 0$ for some agent~$i^*$, the
algorithm presented in Theorem~\ref{th:event-control} is demanding
that events be triggered continuously, i.e., that the control signal
be updated continuously. This happens because the instance one of
agent~$i^*$'s neighbors begins moving, agent~$i^*$ should also be
moving immediately, but the only way to ensure this is to update the
control signal continuously. Since this is not physically possible,
the result of Theorem~\ref{th:event-control} is incomplete until we
can rule out the possibility of Zeno behavior.
}

\begin{remark}\longthmtitle{Zeno behavior and general
    {networked} systems} 
  {\rm The issue pointed out above is not specific to consensus
    problems, and in fact is characteristic of distributed
    event-triggered algorithms operating on networks.  More
    specifically, when a centralized controller is determining when
    events are triggered, this results in a single time schedule for
    which it must be guaranteed that an infinite number of events are
    not triggered in a finite time period. However, when developing a
    distributed event-triggered strategy, individual agents make
    independent decisions regarding when events occur based on partial
    information. This may not only result in many more triggers
    occurring than in the centralized case, but also considerably
    complicates obtaining guarantees about avoiding deadlocks in the
    network. Such analysis usually requires the characterization of
    additional properties of the original algorithm regarding
    robustness to error and the impact of the inter-agent
    interconnections on the evolution of their states.}  \oprocend
\end{remark}

\subsection{Classification of Event-Triggered Consensus
  Algorithms}\label{se:classification}

We have presented above in detail the distributed event-triggered
control problem (Problem~\ref{pr:event-control}) and solution
(Theorem~\ref{th:event-control}).  Since the conception of this
problem and solution, the literature has grown significantly both in
numbers and complexity of the problems and solutions considered. To
help navigate it, our goal here is to identify a number of categories
to systematically classify different problem-solution pairs by their
properties.  For example, we define this particular problem-solution
pair (Problem~\ref{pr:event-control} and
Theorem~\ref{th:event-control}) to have single-integrator dynamics, an
undirected interaction graph, events that trigger control updates,
triggers that are evaluated continuously, and trigger thresholds that
are state-dependent.

In particular, we focus on five main categories to help distinguish
different problem instances and their solutions: \textbf{Dynamics,
  Topology, Trigger Response, Event Detection,} and \textbf{Trigger
  Dependence.} {The first two categories are related to the
  physical problem setup, where Dynamics describes the specific type
  of agent dynamics and Topology captures the type of interactions
  across network agents.  The last three categories are related to the
  capabilities and/or assumptions placed on the agents
  communication/computation abilities.  Trigger Response refers to the
  actions taken by agents in response to an event being triggered,
  Event Detection refers to how events described by triggering
  functions are monitored, and finally, Trigger Dependence refers to
  the arguments and variables that triggering functions depend on.  }
Table~\ref{ta:classifications} summarizes the main distinctions which
are covered in this section in further detail.

We begin by discussing the shortcomings of the
problem-solution pair presented in Section~\ref{se:canonical}, 
and how they can be addressed. The remainder
of this article is then devoted to describing in detail exactly what
the different categories of Table~\ref{ta:classifications} mean, showing exactly how the
different properties change the canonical problem-solution pair
described in Section~\ref{se:canonical}
(Problem~\ref{pr:event-control} and Theorem~\ref{th:event-control}),
and surveying the vast field in terms of these newly proposed
categories and properties. 

%\begin{table*}[h!]
%\begin{center}
%\begin{tabular}{|l|l|l|}
%\hline
%& \textbf{Category} & \textbf{Properties} \\
%\hline
%\dynamics & Dynamics & single-integrator, double-integrator, linear, nonlinear, discrete-time \\
%\hline
%\topology & Topology & undirected, directed balanced, directed rooted spanning tree \\
%\hline
%\triggertype & Trigger Action & control update, broadcast info, request info, edge-based \\
%\hline
%\triggereval & Trigger Evaluation & continuous, periodic, aperiodic \\
%\hline
%\triggerdependence & Trigger Threshold Dependence & state-dependent, time-dependent \\
%\hline
%\inbetween & In-between Triggers & Zero-Order Hold (ZOH), Model Predictive (MP) \\
%\hline
%\uncertainties & Uncertainties & none, state disturbances, sampling noise, unreliable communication \\
%\hline
%\end{tabular}
%\end{center}
%\caption{Different categories for classification of event-triggered consensus problem-solution pairs.}\label{ta:categories}
%\end{table*}

%\begin{figure*}[htb]
%  \centering
%  {\includegraphics[height=.435\linewidth]{scope}}
%  \caption{Different classifications of event-triggered consensus problems detailed in this article.}\label{fig:categories}
%  \vspace*{-1ex}
%\end{figure*}
%\begin{figure*}[htb]
%  \centering
%  {\includegraphics[height=.635\linewidth]{full2}}
%  \caption{Different classifications of event-triggered consensus
%    problems in general, highlighted green boxes correspond to
%    properties this article discusses in
%    detail.}\label{fig:categories2}
%  \vspace*{-1ex}
%\end{figure*}

\section{Event-Triggered Consensus Algorithms}\label{se:main}

In this section we carefully discuss the different types of
event-triggered consensus algorithms outlined in
Section~\ref{se:classification}. We begin by exploring the different
roles a triggering function has on the system. More specifically, we
look at what types of actions agents take in response to a trigger,
how often the triggering functions are evaluated, and what exactly the
triggering functions depend on. {We partition this discussion on
triggers into three categories: Trigger Response, Event Detection, and
Trigger Dependence. 

Let us discuss what capabilities agents physically need to realize
different solutions to Problem~\ref{pr:event-control}.  We note that,
barring the distributed computation aspect, the centralized
event-triggered controller presented in Theorem~\ref{th:central} is a
solution to Problem~\ref{pr:event-control}, where all agent
triggers~$f_i(\cdot)$ are defined as in~\eqref{eq:centraltrigger}.
However, implementing this solution requires all agents to have exact
global state information~$x$ at all times to properly monitor the
function~\eqref{eq:centraltrigger}. Instead,
Theorem~\ref{th:event-control} relaxes this requirement by providing a
local event-triggering function~$f_i(e_i,x_{\NN_i})$ that each
agent~$i$ can monitor with only its neighbors' state
information~$x_{\NN_i}$. However, this solution still requires each
agent~$i$ to have exact state information about their neighbors at all
times. If state information is communicated wirelessly, this means
continuous wireless communication to implement the solution. In the
following, we propose various solutions to
Problem~\ref{pr:event-control} that require less stringent
assumptions.  }

\subsection{\textbf{Trigger Response}}\label{se:triggertype}

\begin{table*}[htb!]
{
  \begin{center}
    \begin{tabular}{|l|l|l|}
      \hline
      \textbf{Category} & \textbf{Properties} & \textbf{Technical Meaning} \\
      \hline
      \begin{tabular}{@{}l@{}l}
        Trigger Response
        \\
        (Section~\ref{se:triggertype})
      \end{tabular} 
      & 
      \begin{tabular}{@{}l@{}l}
        control update only \\ control update and information push  \\
        control update and information pull \\ control update and
        information exchange 
      \end{tabular}
      & 
      \begin{tabular}{@{}l@{}l} 
        $u_i$ updated at event times~$\{t_\ell^i\}_{\ell \in
          \integernonnegative}$ \\ $u_i, \hat{x}_i$ updated at
        event times \\ $u_i, \{\hat{x}_j\}_{j \in \NN_i}$
        updated at event times \\ $u_i, \hat{x}_i, \hat{x}_j$
        (for some $j \in \NN_i)$ updated at event times 
      \end{tabular} 
      \\
      \hline
      \begin{tabular}{@{}l@{}l}
        Event Detection
        \\
        (Section~\ref{se:triggereval})
      \end{tabular} 
      & \begin{tabular}{@{}l@{}l}
        continuous 
        \\
        periodic 
        \\
        aperiodic
      \end{tabular}
      & 
      \begin{tabular}{@{}l@{}l} 
        trigger evaluated at all times~$t
        \in \realnonnegative$ 
        \\
        trigger evaluated periodically~$t \in
        \{0, h_i, 2h_i, \dots \}$
        \\ 
        trigger evaluated
        aperiodically~$t \in \{t_0, t_1, \dots \}$ 
      \end{tabular}
      \\
      \hline
      \begin{tabular}{@{}l@{}l}
        Trigger Dependence
        \\
        (Section~\ref{se:triggerdependence})
      \end{tabular} 
      & 
      \begin{tabular}
        {@{}l@{}l}
        static: state
        \\
        static: time
        \\
        dynamic 
      \end{tabular}
      &
      \begin{tabular}{@{}l@{}l}
        $f_i(\cdot) = f_i(e_i,w_i)$ 
        \\
        $f_i(\cdot) = f_i(e_i,t)$
        \\
        $f_i(\cdot) = f_i(e_i,w_i,\chi_i),
        \quad \dot{\chi}_i = \eta_i(e_i,w_i,\chi_i)$
      \end{tabular} 
      \\
      \hline
      \begin{tabular}{@{}l@{}l}
        Topology
        \\
        (Section~\ref{se:topology})
      \end{tabular} 
      & 
      \begin{tabular}{@{}l@{}l} 
        static \\ dynamic
      \end{tabular}
      & 
      \begin{tabular}{@{}l@{}l}
        communication graph $\commgraph$ is fixed/constant \\
        communication graph~$\commgraph$ is changing over time 
      \end{tabular}
      \\
      \hline
      \begin{tabular}{@{}l@{}l}
      Dynamics
      \\
      (Section~\ref{se:dynamics})
      \end{tabular} 
      &  
      \begin{tabular}{@{}l@{}l}
        single-integrator \\
        double-integrator  \\ linear \\ nonlinear 
      \end{tabular}
      &  
      \begin{tabular}{@{}l@{}} 
        $\dot{x}_i(t) = u_i(t)$ \\
        $\ddot{x}_i(t) = u_i(t)$ \\ $\dot{x}_i(t) = A_i x_i(t) + B_i
        u_i(t)$ \\ $\dot{x}_i(t) = F_i(x_i(t),u_i(t))$
      \end{tabular}
      \\
      \hline
    \end{tabular}
  \end{center}
  \caption{Description of the technical differences between different
    category classifications.}\label{ta:classifications}  
    }
\end{table*}

We begin by discussing the different actions that agents might take in
response to an event being triggered. In the previous section, we
presented event-triggered \emph{control} laws to determine when
control signals should be updated; however, this relies on the
continuous availability of some state information. In particular, each
agent~$i$ requires exact state information about their neighbors~$j
\in \NN_i$ to evaluate the trigger~\eqref{eq:decentraltrigger} and
determine when its control signal $u_i$ should be
updated. {Instead, here we are interested in applying the
  event-triggered paradignm to also drive when communication among
  agents should occur in addition to control updates. We refer to the
  combination of event-triggering for communication and control as
  `event-triggered coordination.'}

% Consequently, we consider 3 new properties to consider within the
% Trigger Actions category~\triggertype: broadcasting (or push-based
% coordination), pull-based coordination, or edge-based
% coordination. We discuss the details of each solution property next.
%

As in Section~\ref{se:canonical}, we assume each agent~$i$ has its own
sequence of event times
$\{ t_\ell^i \}_{\ell \in \integernonnegative}$.  However, these event
times now correspond to when messages are broadcast by agent~$i$; not
just when control signals are updated. At any given time~$t$, let
\begin{align}\label{eq:xhatbroadcast}
  \hat{x}_i(t) = x_i(t_\ell^i) \text{ for } t \in [t_\ell^i, t_{\ell+1}^i )
\end{align} 
be the last \emph{broadcast} state of agent~$i$. Then, at any given
time~$t$, agent $i$ only has access to the last broadcast state
$\hat{x}_j(t)$ of its neighbors $j \in \NN_i$ rather than exact states
$x_j(t)$.

The distributed event-triggered controller is then still given by
\begin{align}\label{eq:decentralcontrol2}
  u_i(t) = - \sum_{j \in \NN_i} (\hat{x}_i(t) - \hat{x}_j(t) ) .
\end{align}
It is important to note here that the latest broadcast state
$\hat{x}_j(t)$ of agent $j \in \NN_i$ appears in the control signal
for agent~$i$ at any time~$t$. This means that when an event is
triggered by a neighboring agent~$j$, agent~$i$ also updates its
control signal accordingly.  As before, let
$e_i(t) = \hat{x}_i(t) - x_i(t)$ be the state measurement error for
agent~$i$. Then, letting $\hat{x} = (\hat{x}_1, \dots, \hat{x}_N)^T$
and $e = (e_1, \dots, e_N)^T$, the closed-loop dynamics of the
controller~\eqref{eq:decentralcontrol2} is again given by
\begin{align}\label{eq:decentraldynamics2}
  \dot{x}(t) = - L \hat{x}(t) = - L ( x(t) + e(t) ).
\end{align}
However, it should be noted that we are now looking for an
event-trigger for each agent~$i$ that does not require exact
information about its neighbors. More specifically, we recall the
result of Theorem~\ref{th:event-control} and notice that the
event-trigger for agent~$i$ depends on the exact state $x_j(t)$ of all
its neighbors $j \in \NN_i$. It was first identified
by~\cite{EK-XW-NH:10} that this solution may not be practical in many
cases, particularly in wireless network settings, as this means agents
must be in constant communication with each other.  {Instead, we are
interested in finding a solution that only depends on the last
broadcast information~$\hat{x}_j(t)$.
%The problem can now be formalized as follows.
%
%\begin{problem}\longthmtitle{Decentralized event-triggered
%    coordination}\label{pr:dec-event}
%  Given a connected, undirected graph~$\commgraph$ and the closed-loop dynamics~\eqref{eq:decentraldynamics2}, find a
%  local event-trigger~$f_i(x_i,\hat{x}_i,\{\hat{x}_j\}_{j \in \NN_i})$
%  for each agent~$i$ such that the sequences of times
%  $\{ t_\ell^i \}_{\ell \in \integernonnegative}$ ensures multi-agent
%  average consensus~\eqref{eq:averageconsensus} is achieved.
%\end{problem}

Specifically, we are now looking for a triggering function~$f_i$ that
only depends on its own state~$x_i$ and the last broadcast state of
its neighbors~$\{\hat{x}_j\}_{j \in \NN_i}$, rather than their true
states~$\{x_j\}_{j \in \NN_i}$.}
Following~\citep{CN-JC:14-acc,EG-YC-HY-PA-DC:13}, to solve this
problem we again consider the Lyapunov function
\begin{align}\label{eq:lyap}
  V(x) = \frac{1}{2}x^T L x.
\end{align}
Given the closed-loop dynamics~\eqref{eq:decentraldynamics2}, we have
\begin{align*}
  \dot{V} = -\TwoNorm{Lx}^2 - x^T L L e,
\end{align*}
just as we did in~\eqref{eq:lyapderivative1}
and~\eqref{eq:lyapderivative2}. However, since we are interested in
identifying conditions for~$\dot{V}$ to be negative in terms of the
most recently broadcast information~$\hat{x}$ instead of actual state
information, we can rewrite this using~$e = \hat{x} - x$ as
\begin{align*}
  \dot{V} = - \TwoNorm{L\hat{x}}^2 + \hat{x}^T L L e. 
\end{align*}

Letting $\hat{z} = L \hat{x} = (\hat{z}_1, \dots, \hat{z}_N)$, it is
easy to see that we can bound
\begin{align*}
  \dot{V} \leq - \left[ \sum_{i=1}^N (1 - a|\NN_i|)\hat{z}_i^2 + \frac{1}{a}
  |\NN_i| e_i^2  \right]
\end{align*}
for all~$a > 0$ following essentially the same steps to arrive
at~\eqref{eq:finalbound}.  Letting~$a \in (0,1/|\NN_i|)$ for all~$i$,
if we can enforce the error of all agents to satisfy
\begin{align}\label{eq:enforcethis}
  e_i^2 \leq \frac{\sigma_i a(1-a|\NN_i|)}{|\NN_i|} \hat{z}_i^2
\end{align}
with $\sigma_i \in (0,1)$ for all times, we have
\begin{align*}
  \dot{V} \leq \sum_{i=1}^N (\sigma_i - 1)(1 - a|\NN_i|) \hat{z}_i^2 
\end{align*}
which is strictly negative for all~$L\hat{x} \neq 0$. {
In order to compute~$\hat{z}_i$, agent~$i$ only needs access to its own and neighbors' broadcast states rather than true states,
\begin{align*}
w_i = \hat{x}_{\NN_i} \triangleq (\hat{x}_i, \{\hat{x}_j\}_{j \in \NN_i} ).
\end{align*} 
}
The following decentralized event-trigger then ensures that~$\dot{V}$ is strictly negative until
consensus is achieved.

\begin{theorem}\longthmtitle{Decentralized event-triggered
    coordination~\citep{EG-YC-HY-PA-DC:13}}\label{th:broadcast1} 
  Given {a connected, undirected graph~$\commgraph$ and} the
  closed-loop dynamics~\eqref{eq:decentraldynamics2}, if event times
  of each agent~$i$ are determined as the times when
  \begin{align}\label{eq:decentraltrigger2}
    f_i(e_i,\hat{x}_{\NN_i} ) \triangleq
    e_i^2 - \frac{\sigma_i a(1-a|\NN_i|)}{|\NN_i|} \hat{z}_i^2 \geq 0,
  \end{align}
  with $0 < a < 1/|\NN_i|$ for all~$i \in \until{N}$, then all
  non-Zeno trajectories of the system asymptotically achieve
  multi-agent average consensus.
\end{theorem}

It is important to note here that unlike the previous
triggers~\eqref{eq:centraltrigger} and~\eqref{eq:decentraltrigger},
agent~$i$'s event is triggered whenever the
inequality~\eqref{eq:decentraltrigger2} is satisfied rather than an
equality. This is because unlike the previous two triggering
functions, this one is discontinuous because it depends on the last
broadcast state~$\hat{x}$ rather than the exact state~$x$, which can
abruptly change anytime an agent triggers an event. In any case, the
point of the trigger is to ensure~\eqref{eq:enforcethis} is satisfied
at all times, which this trigger does (because the instant it is
violated for some agent~$i$, agent~$i$ can immediately set~$e_i = 0$).

Note that Theorem~\ref{th:broadcast1}) and
Theorem~\ref{th:event-control} both solve the exact same problem:
Problem~\ref{pr:event-control}; except we have now considered events
that not only determine when control signals should be updated, but
also when agents should broadcast information to their neighbors. In
more general terms, we refer to this as a \textbf{Trigger Response} of
control signal updates \emph{and} information pushing. Information
pushing, or broadcasting, refers to the action of an agent~$i$ pushing
unsolicited information onto its neighbors~$j \in
\NN_i$. {Table~\ref{ta:classifications} describes the different ways in
which we classify both the problems and the solutions.

Following Table~\ref{ta:classifications}, we can say this
problem-solution pair has single-integrator dynamics, an undirected
interaction graph, events that trigger broadcasts and control updates,
triggers that are evaluated continuously, and trigger thresholds that
are state-dependent. Note, however, that these properties alone are
not enough to uniquely identify a solution to
Problem~\ref{pr:event-control}. Next, we present an alternate solution
to Problem~\ref{pr:event-control} that is described by the exact same
properties listed above.}

Following~\cite{CN-JC:16-auto}, to solve this problem in a different
way we consider a different Lyapunov function,
\begin{align}\label{eq:otherlyap}
  V(x) = \frac{1}{2} (x - \bar{x} \mathbf{1})^T (x - \bar{x}
  \mathbf{1}) ,
\end{align}
where $\bar{x} = \frac{1}{N} \sum_{i=1}^N x_i(0)$ is the average of
all initial conditions.  Then, given the closed-loop
dynamics~\eqref{eq:decentraldynamics2}, we have
\begin{align*}
  \dot{V} = x^T \dot{x} - \bar{x} \mathbf{1}^T \dot{x} = - x^T L
  \hat{x} - \bar{x} \mathbf{1}^T L \hat{x} = - x^T L \hat{x},
\end{align*}
where we have used the fact that the graph is symmetric in the last
equality.  As always, we are interested in finding conditions on the
error~$e$ such that $\dot{V} < 0$ at all times; however, we must now
do it without access to neighboring state information.

As in the previous solution methods, our first step is to
upper-bound~$\dot{V}$ to find conditions to ensure it is never
positive. The following result from~\citep{CN-JC:16-auto} can then be
used to find conditions that ensure~$\dot{V} < 0$ at all times.

\begin{lemma}[\citep{CN-JC:16-auto}]\label{le:bound}
  Given~$V(x) =
  \frac{1}{2}(x-\bar{x}\mathbf{1})^T(x-\bar{x}\mathbf{1})$
  and the closed-loop dynamics~\eqref{eq:decentraldynamics2},
  \begin{align*}
    \dot{V} \leq \sum_{i=1}^N e_i^2 |\NN_i| - \sum_{j \in \NN_i}
    \left( \frac{1}{4}(\hat{x}_i-\hat{x}_j)^2 \right).  
  \end{align*} 
\end{lemma}
% \begin{pf}
%   In the appendix.
% \end{pf}

Leveraging Lemma~\ref{le:bound}, it is easy to see that if we can
enforce the error of all agents to satisfy
\begin{align}\label{eq:enforcethisalways}
  e_i^2 \leq \sigma_i \frac{1}{4 | \NN_i |} \sum_{j \in \NN_i}
  (\hat{x}_i - \hat{x}_j)^2
\end{align}
with $\sigma_i \in (0,1)$ for all times, we have
\begin{align}\label{eq:thenthis}
  \dot{V} \leq \sum_{i=1}^N \frac{\sigma_i - 1}{4} \sum_{j \in \NN_i}
  (\hat{x}_i - \hat{x}_j)^2,
\end{align}
which is strictly negative for all $L\hat{x} \neq 0$. For simplicity,
we use the shorthand notation
\begin{align*}
  \hat{\phi}_i = \sum_{j \in \NN_i} (\hat{x}_i - \hat{x}_j)^2 
\end{align*}
in the definition of the following decentralized event-trigger, which
ensures that this is satisfied at all times.

\begin{theorem}\longthmtitle{Decentralized event-triggered
    coordination~\citep{CN-JC:16-auto}}\label{th:broadcast2}
  Given {a connected, undirected graph~$\commgraph$ and} the
  closed-loop dynamics~\eqref{eq:decentraldynamics2}, if the event
  times of each agent~$i$ are determined as the times when
  \begin{align}
    f_i(e_i,\hat{x}_{\NN_i}) \triangleq e_i^2 -
    \sigma_i \frac{1}{4 |\NN_i|} \hat{\phi}_i \geq 0,
  \end{align}
  for all~$i \in \until{N}$, then all non-Zeno trajectories of the
  system achieve multi-agent average consensus.
\end{theorem}

%\begin{table}[htb]
%  \centering
%  \framebox[.9\linewidth]{\parbox{.85\linewidth}{%
%      \parbox{\linewidth}{At times $t \in [t_\ell^i,t_{\ell+1}^i)$,
%        agent~$i$ (continuously) performs:} 
%      \vspace*{-1.5ex}
%      \begin{algorithmic}[1]
%      %\STATE set $z_i(t) = \sum_{j \in \NN_i} (x_i(t) - x_j(t))$
%      \STATE set $\hat{x}_i(t) = x_i(t_\ell^i)$
%      \STATE set $e_i(t) = \hat{x}_i(t) - x_i(t)$
%        \IF{$e_i(t)^2 \geq \sigma_i \frac{1}{4 | \NN_i |} \sum_{j \in \NN_i} (\hat{x}_i(t) - \hat{x}_j(t))^2$}  
%        \STATE set $t_{\ell+1}^i = t$
%        \STATE broadcast $\hat{x}_i(t) = x_i(t_{\ell+1}^i)$ to neighbors $j \in \NN_i$
%        \STATE set $\ell = \ell + 1$
%        \ENDIF
%      \STATE set $u_i(t) = - \sum_{j \in \NN_i} (\hat{x}_i(t) - \hat{x}_j(t))$
%      \end{algorithmic}}}
%  \caption{Decentralized event-triggered coordination (state-dependent).}\label{tab:algorithm5}
%\end{table}

Note that both solutions presented in Theorem~\ref{th:broadcast1} and
Theorem~\ref{th:broadcast2} are classified the same way according to
the categories in Table~\ref{ta:classifications}, because they can be
implemented under the same assumptions on the agent capabilities. In
this case agents need to be able to receive broadcasted information
from their neighbors on some connected, undirected communication
graph. \changes{Note that both of these solutions are still
  technically incomplete as they do not rule out the possibility of
  Zeno behavior. In fact, the original design in~\cite{CN-JC:16-auto}
  includes an additional trigger to ensure that Zeno behavior does not
  occur.}

We have now discussed two different actions an agent might take in
response to an event being triggered: control updates and broadcasting
a message. Instead, one could imagine other types of actions resulting
from an event being triggered as well. In particular, in
Table~\ref{ta:classifications} we highlight the possibility of events
triggering \emph{information pulls} or \emph{exchanges}, rather than
broadcasts (or information pushes).  More specifically, we refer to an
information push by an agent as a broadcast message that can be
received by all its neighbors. Instead, we can think of an information
pull by an agent as a request for updated information to its
neighbors. That is, when an event is triggered by an agent, rather
than telling its neighbors its current state, it instead requests
state information from its neighbors.  {This idea is related to
self-triggered control design~\citep{AA-PT:10,WPMHH-KHJ-PT:12}, where
the decision maker, at each event time, immediately schedules its next
event time with the information available at the current event time
(rather than continuously monitoring a triggering condition as one
normally does in event-triggered control). In self-triggered
coordination with wireless communication, agents use their current
information to determine when in the future they need to acquire new
information from
others~\citep{CN-JC:11-auto,CDP-PF:13,YF-LL-GF-YW:15,EH-DEQ-EGWP-HS-KHJ:15}.}

In addition to information pushes or pulls, one could also imagine
scenarios in which an information exchange or swap may be more
practical. Applying event-triggered ideas to gossiping protocols has
recently been called edge-based event-triggered
coordination~\citep{CDP-PF:13,CDP-RP:17,MC-FX-LW:15,WX-DWCH-JZ-BC:16,BW-FX-MD:17,GD-FX-LW:17}. In
this setting events are triggered along specific edges of a graph
rather than its nodes. More specifically, rather than a single agent
sending information to or requesting information from \emph{all} its
neighbors at once, events are instead triggered at the link level
which drives a direct agent-to-agent information exchange.

\subsection{\textbf{Event Detection}}\label{se:triggereval}

Having discussed the different types of actions that an event can
drive, until now we have assumed that the triggering functions can be
monitored continuously.  This is troublesome since this is technically
not possible for cyber-physical systems.  Unfortunately, even if
continuous event detection were possible, most of the algorithms
presented in the article thus far are not guaranteed to avoid Zeno
behaviors making them risky to implement on real systems. In fact, as
mentioned above, until Zeno behavior is guaranteed not to
occur in the system, the convergence results of
Theorems~\ref{th:event-control},~\ref{th:broadcast1},
and~\ref{th:broadcast2} are not even valid.

Consequently, being able to properly rule out the existence of Zeno
behavior in an event-triggered consensus problem is both subtle and
critical for its correctness. Recalling Definition~\ref{def:zeno},
Zeno behavior is defined as having an infinite number of events
triggered in any finite time period. Unfortunately, it turns out that
in all the algorithms presented so far, it is not guaranteed that
Zeno behavior will not occur.

Note that in some cases the algorithms can be slightly modified to
theoretically avoid Zeno behavior, but even in these cases it turns
out that the time between two events generated by a single agent may
be arbitrarily small, see e.g.,~\cite{CN-JC:16-auto}. {More
  specifically, even if it can be guaranteed that an infinite number
  of events are not triggered in any finite time period, the time
  between two events might not have a uniform lower bound. This means
  that even with a non-Zeno guarantee, this is still troublesome from
  an implementation viewpoint because an agent's hardware/software
  physically cannot keep up with how quickly events are being
  generated, cf. Section~\ref{se:zeno}.  }

Motivated by this discussion, researchers have considered enforcing a
minimum time between events as a more practical constraint for
event-triggered solutions.  Until now we have assumed that all
event-triggers can be evaluated continuously. Or more specifically we
say that the \textbf{Event Detection} occurs continuously.  That is,
the exact moment at which a triggering condition is met, an action
(e.g., state broadcast and control signal update) is carried
out. However, as mentioned above this may be an unrealistic assumption
when considering actual digital implementations. More specifically, a
real device cannot continuously evaluate whether a triggering
condition has occurred or not.

This observation motivates the need for relaxing the continuous event
detection requirement and instead determine a discrete set of times at
which triggering functions should be evaluated.  The most natural way
to approach this is to study sampled-data (or periodically checked)
event-triggered coordination strategies.

Specifically, given a sampling period $h \in \realpositive$, we let
$\{ t_{\ell'} \}_{\ell' \in \integernonnegative}$, where
$t_{\ell'+1} = t_\ell' + h$, denote the sequence of times at which
agents evaluate decisions about whether to broadcast their state to
their neighbors or not. This type of design is more in line with the
constraints imposed by real-time implementations, where individual
components work at some fixed frequency, rather than continuously. An
inherent and convenient feature of this strategy is the automatic lack
of Zeno behavior (since inter-event times are naturally lower bounded
by~$h$).

Under the new framework we still have familiar equations. In
particular, the control law of each agent is still given by
\begin{align}\label{eq:periodiccontrol}
  u_i(t) = - \sum_{j \in \NN_i} (\hat{x}_i - \hat{x}_j),
\end{align}
which means we still have the same closed-loop dynamics,
\begin{align}\label{eq:periodicdynamics}
  \dot{x}(t) = -L \hat{x}(t) = - L(x(t) + e(t)).
\end{align}
The difference only shows up when considering when broadcasts
occur. That is, we now have
\begin{align}\label{eq:xhatperiodic}
  \hat{x}_i(t) = x_i(t_\ell^i) \text{ for } t \in [t_\ell^i, t_{\ell+1}^i),
\end{align}
just as we did in~\eqref{eq:xhatbroadcast} when considering continuous
event detection, except now the event
times~$\{t_\ell^i\} \subset \{0, h, 2h, \dots \}$ can only occur at
discrete time-instances.
%
%\begin{problem}\longthmtitle{Decentralized periodic event-triggered
%    coordination}\label{pr:periodic}
%  Given a connected, undirected graph~$\commgraph$ and the closed-loop dynamics~\eqref{eq:periodicdynamics}, find a
%  local event-trigger~$f_i(x_i,\hat{x}_i,\{\hat{x}_j\}_{j \in \NN_i})$
%  for each agent~$i$ such that the sequences of times
%  $\{ t_\ell^i \}_{\ell \in \integernonnegative} \subset \{0, h, 2h,
%  \dots \}$
%  ensures multi-agent average consensus~\eqref{eq:averageconsensus} is
%  achieved.
%\end{problem}

{To again solve Problem~\ref{pr:event-control} by now relaxing
  the continuous monitoring requirement, let us begin by revisiting
  the result of Theorem~\ref{th:broadcast2}, where we used the
  Lyapunov function~\eqref{eq:otherlyap},}
\begin{align}\label{eq:otherlyapagain}
  V(x) = \frac{1}{2} (x - \bar{x} \mathbf{1})^T (x - \bar{x}
  \mathbf{1}) ,
\end{align}
and Lemma~\ref{le:bound} to find conditions on the error~$e_i$ such
that $\dot{V}$ was always negative. Since the Lyapunov function here
and closed-loop dynamics~\eqref{eq:periodicdynamics} has the exact
same functional form as before, Lemma~\ref{le:bound} still holds as
well. That is, given the Lyapunov function~\eqref{eq:otherlyapagain}
and closed-loop dynamics~\eqref{eq:periodicdynamics}, we have the
upper-bound
\begin{align}\label{eq:enforcealways}
  \dot{V} \leq \sum_{i=1}^N e_i^2 |\NN_i| - \sum_{j \in \NN_i} \left(
  \frac{1}{4} (\hat{x}_i - \hat{x}_j)^2 \right) . 
\end{align}
Then, just as before, we want to find conditions on~$e_i$ such that
this is always negative.  However, the issue now is that we can only
generate events at the discrete times~$\{0, h, 2h, \dots\}$, meaning
we can only reset the error~$e_i(t) = 0$ at these specific times
rather than any time~$t \in \realnonnegative$. Thus, we must now find
a triggering condition that is only evaluated at these sampling times,
but still guarantees that the righthand side
of~\eqref{eq:enforcealways} is negative \emph{for all $t \in
  \realnonnegative$}.

Intuitively, as long as the sampling period~$h$ is small enough, the
closed-loop system with a periodically checked event-triggering
condition will behave similarly to the system with triggers being
evaluated continuously. Interestingly, two different groups have
developed two different, albeit similar, algorithms based on the same
Lyapunov function~\eqref{eq:otherlyapagain} using two different ways
of upper-bounding its time-derivative.  We omit the details here but
present the two solutions to this problem based on the
works~\citep{XM-TC:13,CN-JC:16-auto} next.

\begin{theorem}\longthmtitle{Periodic event-triggered
    coordination~\citep{XM-TC:13}}\label{th:sampled-data1}  
  Given {a connected, undirected graph~$\commgraph$ and} the
  closed-loop dynamics~\eqref{eq:periodicdynamics}, if the event times
  of each agent~$i$ are determined as the times $t' \in \{0, h, 2h,
  \dots \}$ when
  \begin{align*}
    f_i(e_i, \hat{x}_{\NN_i}) \triangleq e_i^2 - \sigma_i
    \hat{z}_i^2\geq 0,
  \end{align*}
  and $h \in \realpositive$ and $\sigma_\text{max}$ satisfy
  \begin{align}\label{eq:bound}
    h \leq \frac{1}{2 \lambda_N} \text{ and } \sigma_\text{max} <
    \frac{1}{\lambda_N^2} ,
  \end{align}
  where $\sigma_\text{max} = \max_{i \in \until{N}} \sigma_i$, then the
  system achieves multi-agent average consensus.
\end{theorem}

\begin{theorem}\longthmtitle{Periodic event-triggered
    coordination~\citep{CN-JC:16-auto}}\label{th:sampled-data2} 
  Given {a connected, undirected graph~$\commgraph$ and} the
  closed-loop dynamics~\eqref{eq:periodicdynamics}, if the event times
  of each agent~$i$ are determined as the times $t' \in \{0, h, 2h,
  \dots \}$ when
  \begin{align}\label{eq:trigger4}
    f_i(e_i, \hat{x}_{ \NN_i}) \triangleq e_i^2 - \sigma_i
    \frac{1}{4 |\NN_i|} \hat{\phi}_i \geq 0,
  \end{align}
  and $h \in \realpositive$ and $\sigma_\text{max}$ satisfy
  \begin{align}\label{eq:bound2}
    % \sigma_\text{max} + 2hC < 1,
    \sigma_{\max} + 4 h |\NN_\text{max} |^2 < 1,
  \end{align}
  where $|\NN_\text{max} | = \max_{i \in \until{N}} |\NN_i |$,
  then the system achieves multi-agent average consensus.
\end{theorem}

Note that the results of Theorems~\ref{th:sampled-data1}
and~\ref{th:sampled-data2} guarantee that \emph{all trajectories} of
the systems can achieve multi-agent average consensus under their
respective conditions since Zeno executions are trivially ruled out
because agents can only trigger an event at most every~$h > 0$
seconds. \changes{This trivially gives us a positive MIET~$\tau^\text{min} = h$.}	
We also note here the small difference in the triggering
functions and conditions on~$h$ and~$\sigma_\text{max}$ for
convergence of these results are a result of different ways of
upper-bounding~$\dot{V}$ and ultimately being able to guarantee that
$\dot{V} < 0$ at all times. In particular, we note that the conditions
for guaranteeing convergence in Theorem~\ref{th:sampled-data1} are
less conservative but requires the algebraic connectivity of the
communication graph; whereas the conditions for guaranteeing
convergence in Theorem~\ref{th:sampled-data2} may be more strict but
easier to compute by the agents themselves.

A drawback of these and similar
solutions~\citep{XM-LX-YCS-CN-GJP:15,YF-YY-YZ:16} is that the
period~$h$ must be the same for all agents, requiring synchronous
action. This assumption may be restrictive in practical scenarios
where data cannot be consistently acquired.  Instead, it seems
desirable to develop asynchronous versions of these solutions or, more
generally, solutions where the \textbf{Event Detection} occurs
aperiodically rather than continuously or periodically. We are only
aware of a few recent works that have begun investigating the
asynchronism issue~\citep{XM-LX-YCS:17,YL-CN-ZT-QL:17,GD-FX-LW:17,YL-CN-ZT-QL:17}.

More specifically, new algorithms may be required to consider the case
of aperiodic sampled-data event detection, or even self-triggered
event detection. In the former case agents would obtain samples at
different instances of time, and then take appropriate actions in
response. In the latter case, one could imagine a scenario where
agents are not only responsible for determining when communication
should occur, but also when local samples should be taken. In this
case it may be useful to consider self-triggered sampling combined
with event-triggered communication and control. More specifically, the
agents would determine by themselves when future samples should be
taken, and then event-decisions should be made based on the taken
samples.

\subsection{\textbf{Trigger Dependence}}\label{se:triggerdependence}

{We have now discussed what agents should do in response to a trigger
and how carefully these triggers need to be monitored. We are now
interested in studying what these triggering functions should actually
depend on and why. In particular, we have only considered triggering
functions so far that depend on locally available information and no
exogenous signals. In this section we present the difference between
static and dynamic triggering functions.  A static triggering function
means that the trigger only depends on currently available information
(i.e., memoryless), whereas a dynamic triggering function may depend
on additional internal dynamic variables.

We begin with static time-dependent triggering functions, rather than
the state-dependent ones we have used until now.} Let us now revisit
Problem~\ref{pr:event-control} again, except this time we are
interested in designing an event-trigger threshold that is
time-dependent rather than state-dependent. The time-dependent
event-trigger to solve this problem was first developed
by~\cite{GSS-DVD-KHJ:13} and is presented next.

\begin{theorem}\longthmtitle{Decentralized event-triggered
    coordination (time-dependent)~\citep{GSS-DVD-KHJ:13}}\label{th:time}
  Given {a connected, undirected graph~$\commgraph$ and} the
  closed-loop dynamics~\eqref{eq:decentraldynamics2}, if the event
  times of each agent~$i$ are determined as the times when
  \begin{align}
    f_i(e_i,t) \triangleq \TwoNorm{e_i} - (c_0 + c_1 e^{-\alpha
    t} ) = 0,
  \end{align}
  with constants $c_0, c_1 \geq 0$ and $c_0 + c_1 > 0$, then all
  non-Zeno trajectories of the system reach a neighborhood of
  multi-agent average consensus upper-bounded by
  \begin{align*}
    r = \TwoNorm{L} \sqrt{N} c_0 / \lambda_2(L) .
  \end{align*}
  Moreover, if $c_0 > 0$ or $0 < \alpha < \lambda_2(L)$, then the closed-loop
  system does not exhibit Zeno behavior.
\end{theorem}
%\begin{pf}
%  In the appendix. \qed
%\end{pf}

This solution uses a triggering function whose threshold depends on
time rather than state. Thus, we can say this problem-solution pair
has single-integrator dynamics, an undirected interaction graph,
events that trigger broadcasts and control updates, triggers that are
evaluated continuously, and trigger thresholds that are time-dependent
rather than state-dependent. 

{We point out here the closely related notion of
  ``event-triggered mechanism'' (ETM), as presented
  in~\citep{DPB-WPMHH:14}, where three classes are presented:
  relative, absolute, and mixed. Here, we have proposed slightly more
  general classes of trigger dependencies such that the relative ETM
  is a special case of our state-dependent triggers. Similarly, the
  absolute ETM (or constant triggering threshold) is a special case of
  the time-dependent trigger of Theorem~\ref{th:time} with $c_1 = 0$,
  which were among the first types of event-trigger thresholds
  considered in network settings unrelated to
  consensus~\citep{MM:06,MZ-CGC:10}. The mixed ETM is a combination of
  these two triggers, but we do not discuss the distinction in this
  article.

  Constant thresholds (or absolute ETMs) give two main advantages. The
  first is their simplicity to implement, and the second is that it is
  generally easy to rule out the possibility of Zeno behavior for
  them. Since the threshold is a constant, it usually takes some
  nonzero minimum amount of time for the error to be able to reach
  this threshold from zero. Note that this is evident in the result of
  Theorem~\ref{th:time} which guarantees for~$c_0 > 0$ that Zeno
  behaviors do not occur. The drawback is that the constant thresholds
  generally do not generate events at times that align well with the
  evolution of the task at hand, and hence, the price we pay is that
  one is not able to guarantee exact convergence all the way to the
  desired states.  } This is discussed in detail next as a particular
case of the time-dependent algorithm presented above.

%{\color{green}
%\begin{table}[htb]
%  \centering
%  \framebox[.9\linewidth]{\parbox{.85\linewidth}{%
%      \parbox{\linewidth}{At times $t \in [t_\ell^i,t_{\ell+1}^i)$, agent~$i$ (continuously) performs:}
%      \vspace*{-1.5ex}
%      \begin{algorithmic}[1]
%      %\STATE set $z_i(t) = \sum_{j \in \NN_i} (x_i(t) - x_j(t))$
%      \STATE set $\hat{x}_i(t) = x_i(t_\ell^i)$
%      \STATE set $e_i(t) = \hat{x}_i(t) - x_i(t)$
%        \IF{$|e_i(t)| = c_0 + c_1 e^{- \alpha t}$}  
%        \STATE set $t_{\ell+1}^i = t$
%        \STATE broadcast $\hat{x}_i(t) = x_i(t_{\ell+1}^i)$ to neighbors $j \in \NN_i$
%        \STATE set $\ell = \ell + 1$
%        \ENDIF
%      \STATE set $u_i(t) = - \sum_{j \in \NN_i} (\hat{x}_i(t) - \hat{x}_j(t))$
%      \end{algorithmic}}}
%  \caption{Decentralized event-triggered coordination (time-dependent).}\label{tab:algorithm4}
%\end{table}
%}

Both the advantages and disadvantages of the event-triggered
coordination law with the time-dependent triggers
proposed in Theorem~\ref{th:time} come from the tunable design
parameters $c_0, c_1,$ and $\alpha$, which play  important roles
in the performance of the algorithm (e.g., convergence speed and
amount of events triggered). For example, setting $c_1$ small and
$\alpha$ large increases the convergence rate at the cost of more
events being triggered, whereas setting a large $c_0$ reduces the
number of events being triggered at the cost of not being able to
converge exactly to the initial average. These parameters can then be
tuned to give a desired balance between performance and efficiency.
Another advantage of the time-dependent triggers are their simplicity
to design and implement. 

Unfortunately, there are also some physical limits to how these
parameters can be tuned to guarantee Zeno behaviors do not occur. For
example, if $\alpha$ is set too high we may be asking the system to
converge faster than is physically possible, leading to an infinite
number of events being generated in a finite time. In particular we
focus our discussion here on the parameters $c_0$ and $\alpha$ and
their effects on convergence and possible Zeno behaviors. We begin
with the more desirable $c_0 = 0$ case, as in this case the result of
Theorem~\ref{th:time} states that the system will asymptotically
achieve exact multi-agent average consensus as defined
in~\eqref{eq:averageconsensus}. However, in this case we require
$\alpha < \lambda_2(L)$ to guarantee Zeno behaviors can be avoided
and, unfortunately, $\lambda_2(L)$ is a global quantity that requires
knowledge about the entire communication topology to compute. {There
  are indeed methods for estimating this quantity in a distributed way
  (see e.g.,~\cite{RA-GS-DVD-CS-KHJ:12,PY-RAF-GJG-KML-SSS-RS:10}) but
  we do not discuss this here.}  On the other hand, when $c_0 > 0$ we
can guarantee that Zeno behaviors are avoided regardless of our choice
of $\alpha$; however, we lose the exact asymptotic convergence
guarantee. Note that in the case of constant triggers, i.e.,
$c_1 = 0$, we must have~$c_0 > 0$. That means in these cases we can
only guarantee convergence to a neighborhood of the desired average
consensus state rather than exact average consensus.

As a result of the above discussion, we see that it is difficult for
the agents to choose the parameters $c_0, c_1,$ and $\alpha$ without
global knowledge to ensure asymptotic convergence to the average
consensus state while also guaranteeing Zeno executions are
avoided. On the other hand, using state-dependent triggers might be
more risky to implement as it is generally harder to rule out the
possibility of Zeno behaviors.

{Referring back to Table~\ref{ta:classifications}, we have now
  discussed both types of static triggering functions.
  % More specifically, we have only discussed events that are defined
  % by a single static function.
  Some works have also considered hybrid or mixed event-time driven
  coordination, where events may be generated by both state and time
  events~\citep{FX-TC:16,FX-TC-HG:16,DPB-WPMHH:14,ZS-NH-BDOA-ZD:16}. Beyond events
  generated as functions of time or state, other works have also
  considered input-based events that depend on the control signal
  being used~\citep{YW-XM-LX-RL-HS-ZW:16,AD-DL-DVD-KHJ:16}. While the
  analysis is slightly different, the intuitive idea is similar. When
  considering state-based events, we generally trigger an event when
  the error between the true state and state currently used by the
  controller exceeds some threshold. However, ultimately what matters
  in a control system is the signal being used. Thus, these
  input-based triggering functions instead define an input error
  between the actual input being used and the desired input if exact
  state information was available. Thus, even if the state error is
  large, these algorithms do not trigger a controller update until the
  input error exceeds some threshold.

  Beyond static triggers, the idea of a dynamic event-triggering
  function has recently been applied as a promising method to rule out
  Zeno
  behavior~\citep{VSD-WPMHH:15,AG:17,XY-KL-DVD-KHJ:17,VSD-DPB-WPMHH:17}.
  In this case the triggering function~$f(\cdot)$ depends on an
  additional, internal dynamic variable with its own dynamics that can
  be designed separately.
  
  Let us revisit Problem~\ref{pr:event-control}, except that now we
  aim to design a dynamic event-triggered coordination strategy that
  can guarantee average consensus from all initial conditions with no
  global information and including non-Zeno guarantees. This result is
  presented next.

  \begin{theorem}\longthmtitle{Dynamic event-triggered
      control~\citep{XY-KL-DVD-KHJ:17}}\label{th:dynamic}
    Given {a connected, undirected graph~$\commgraph$ and} the
    closed-loop dynamics~\eqref{eq:decentraldynamics}, if the event
    times of each agent~$i$ are determined by
    \begin{subequations}\label{eq:dynamictrigger}
      \begin{align}
        f_i(e_i, \hat{x}_{ \NN_i}, \chi_i ) &\triangleq
        |\NN_i| e_i^2 - \frac{ \sigma_i}{4} \hat{\phi}_i - \chi_i \geq
        0,
        \\
        \changes{\dot{\chi}_i(e_i, \hat{x}_{\NN_i}, \chi_i)} &= -\chi_i + \frac{\sigma_i}{4} \hat{\phi}_i - e_i^2   ,
      \end{align}      
    \end{subequations}
    with $\chi_i(0) > 0$ and $\sigma_i
    \in [0,1)$ for all~$i \in \until{N}$, then the system
    asymptotically achieves multi-agent average consensus.
  \end{theorem}
}

\begin{table*}[tbh!]
  {
  \begin{center}
    \begin{tabular}{|l|l|l|}
      \hline
      \textbf{Solution classification} & \textbf{Triggering mechanism}
      & \textbf{Properties} 
      \\
      \hline 
      \begin{tabular}{@{}l@{}}
        Theorem~\ref{th:event-control}~\citep{DVD-EF-KHJ:12} 
        \\ 
        \textbf{Trigger Response}: control updates
        \\
        \textbf{Event
          Detection}: continuous 
        \\
        \textbf{Trigger Dependence}:
        static: state 
      \end{tabular}
      & 
      \begin{tabular}{@{}l@{}} 
        $f_i(e_i, x_{\NN_i} ) \triangleq e_i^2 - \frac{
          \sigma_i a (1 - a|\NN_i|)}{|\NN_i|} z_i^2 = 0$ 
      \end{tabular}
      & 
      \begin{tabular}{@{}l@{}} 
        requires continuous monitoring 
        \\ 
        \;of neighbors~$x_{\NN_i}$;
        \\
        no non-Zeno guarantee 
      \end{tabular}
      \\  
      \hline
      \begin{tabular}{@{}l@{}} 
        Theorem~\ref{th:broadcast1} \citep{EG-YC-HY-PA-DC:13} 
        \\
        Theorem~\ref{th:broadcast2} \citep{CN-JC:16-auto}
        \\
        \textbf{TR}: control updates, info push
        \\
        \textbf{ED}:
        continuous
        \\ 
        \textbf{TD}: static: state
      \end{tabular} 
      &
      \begin{tabular}{@{}l@{}} 
        $f_i(e_i, \hat{x}_{\NN_i} )
        \triangleq e_i^2 - \frac{\sigma_i a(1-a|\NN_i|)}{|\NN_i|}
        \hat{z}_i^2 \geq 0$
        \\
        \hspace*{11ex} or
        \\ 
        $f_i(e_i,\hat{x}_{\NN_i} ) \triangleq e_i^2 - \sigma_i \frac{1}{4
          |\NN_i|} \hat{\phi}_i \geq 0$
      \end{tabular}
      & \begin{tabular}{@{}l@{}} 
        no non-Zeno guarantee
      \end{tabular}
      \\
      \hline
      \begin{tabular}{@{}l@{}} 
        Theorem~\ref{th:sampled-data1} \citep{XM-TC:13}
        \\
        Theorem~\ref{th:sampled-data2} \citep{CN-JC:16-auto}
        \\ 
        \textbf{TR}: control updates, info push
        \\
        \textbf{ED}: periodic
        \\ 
        \textbf{TD}: static: state
      \end{tabular}
      &
      \begin{tabular}{@{}l@{}}
        (Only at times~$t \in \{0, h, 2h, \dots \}$)
        \\
        $f_i(e_i, \hat{x}_{\NN_i}) \triangleq e_i^2 - \sigma_i 
        \hat{z}_i^2\geq 0$ 
        \\ 
        \hspace*{11ex} or 
        \\
        $ f_i(e_i,\hat{x}_{\NN_i}) \triangleq e_i^2 - \sigma_i
        \frac{1}{4 |\NN_i|} \hat{\phi}_i \geq 0$  
      \end{tabular}
      & 
      \begin{tabular}{@{}l@{}} 
        \changes{positive MIET~$\tau^\text{min} = h$;}
        \\
        requires
        synchronous period~$h > 0$
        \\
        \; to guarantee convergence 
      \end{tabular}
      \\
      \hline
      \begin{tabular}{@{}l@{}} 
        Theorem~\ref{th:time} \citep{GSS-DVD-KHJ:13}
        \\
        \textbf{TR}: control updates, info push
        \\ \textbf{ED}: continuous
        \\ \textbf{TD}: static: time 
      \end{tabular}
      &
      \begin{tabular}{@{}l@{}} 
        $f_i(e_i,t) \triangleq \TwoNorm{e_i} - (c_0 + c_1 e^{-\alpha
          t} ) = 0$
      \end{tabular}
      & 
      \begin{tabular}{@{}l@{}} 
        requires
        algebraic connectivity~$\lambda_2$
        \\ 
        \; to guarantee non-Zeno ;\\
        \changes{no positive MIET}
      \end{tabular}
      \\
      \hline
      \begin{tabular}{@{}l@{}} 
        Theorem~\ref{th:dynamic}
        \citep{XY-KL-DVD-KHJ:17}
        \\
        \textbf{TR}: control updates, info
        push
        \\ 
        \textbf{ED}: continuous
        \\ 
        \textbf{TD}: dynamic 
      \end{tabular}
      &
      \begin{tabular}{@{}l@{}}
        $f_i(e_i,\hat{x}_{\NN_i}, \chi_i ) \triangleq |\NN_i| e_i^2 -
        \frac{\sigma_i}{4}  \hat{\phi}_i  - \chi_i \geq 0$ 
        \\ 
        \hspace*{3ex} $\dot{\chi}_i = -\chi_i + 
          \frac{\sigma_i}{4}  \hat{\phi}_i - e_i^2$ 
      \end{tabular}
      & 
      \begin{tabular}{@{}l@{}} 
        guarantees
        non-Zeno; \\ \changes{no positive MIET}
      \end{tabular}
      \\
      \hline
      \begin{tabular}{@{}l@{}} 
        \changes{Theorem~\ref{th:dynamic2} \\
        \citep{JB-CN:19-acc}}
        \\
        \changes{\textbf{TR}: control updates, info
        push}
        \\ 
        \changes{\textbf{ED}: continuous}
        \\ 
        \changes{\textbf{TD}: dynamic }
      \end{tabular}
      &
      \begin{tabular}{@{}l@{}}
        \changes{$f_i(\chi_i) \triangleq -\chi_i \geq 0$} \\ 
        \hspace*{3ex} \changes{$\dot{\chi}_i = \min \left\lbrace -1 , \frac{\hat{\phi}_i}{e_i^2} + 2(\chi_i + 1) \frac{\hat{z}_i}{e_i} - 1 \right\rbrace$}
      \end{tabular}
      & 
      \begin{tabular}{@{}l@{}} 
        \changes{positive MIET~$\tau_i^\text{min} =$ }\\
        $\frac{1}{\sqrt{|\NN_i|}} \left[ \operatorname{atan}(2 \sqrt{|\NN_i|}) - \operatorname{atan}(\sqrt{|\NN_i|}) \right] $ 
      \end{tabular}
      \\
      \hline
    \end{tabular} 
    \\
  \end{center}
  \caption{Summary of solutions to the decentralized event-triggered  
      consensus problem, cf. Problem~\ref{pr:event-control}, discused
      in this article. Note that 
    all these solutions assume that the communication
    \textbf{Topology} is undirected and connected and the
    \textbf{Dynamics} of each agent are single-integrators. 
    % For convenience, we let $x_{\NN_i} = \{x_i, \{ x_j \}_{j \in
    % \NN_i} \}$ and $\hat{x}_{\NN_i} = \{ \hat{x}_i, \{ \hat{x}_j \}_{j
  % \in \NN_i} \}$ denote the exact and last broadcast state
  % information of agent~$i$ and its neighbors~$j \in \NN_i$. Recall
  % that $z_i = (Lx)_i = \sum_{j \in \NN_i} (x_i - x_j)$, $\hat{z}_i =
  % (L\hat{x})_i = \sum_{j \in \NN_i} (\hat{x}_i - \hat{x}_j)$, and
  % $\hat{\phi}_i = \sum_{j \in \NN_i} (\hat{x}_i - \hat{x}_j)^2.$}
    Table~\ref{ta:taxonomy} recalls all the relevant terms.}\label{ta:newdetails} 
}
\end{table*}

Theorem~\ref{th:dynamic} fully solves Problem~\ref{pr:event-control}
in a distributed way. Notably, this solution does not require agents
to have any global information to implement the algorithm, and
guarantees convergence to the desired consensus state by also
guaranteeing Zeno behavior does not occur along any
trajectory. \changes{However, it should be noted that while this
  solution theoretically solves Problem~\ref{pr:event-control}, it
  does not guarantee the existence of a positive
  MIET~$\tau^\text{min}$, which poses problems for practical
  implementation, as discussed in Section~\ref{se:zeno}. More
  recently,~\cite{JB-CN:19-acc} have developed a new dynamic
  triggering strategy that guarantees a positive MIET for each agent,
  providing a complete and implementable solution to
  Problem~\ref{pr:event-control}. This result is formalized next.

  \begin{theorem}\longthmtitle{Dynamic event-triggered control with a
      positive MIET~\citep{JB-CN:19-acc}}\label{th:dynamic2}
    Given a connected, undirected graph~$\commgraph$ and the
    closed-loop dynamics~\eqref{eq:decentraldynamics}, if the event
    times of each agent~$i$ are determined by
    \begin{align}
      f_i(\chi_i) &\triangleq -\chi_i \geq 0, \\
      \dot{\chi}_i(e_i,\hat{x}_{\NN_i},\chi_i) &=
      \min \left\lbrace -1, \frac{\hat{\phi}_i}{e_i^2} + 2(\chi_i + 1) \frac{\hat{z}_i}{e_i} -1 \right\rbrace,
    \end{align}
    with~$\chi_i(t_\ell^i) \triangleq 1$ for all~$\ell \in
    \integernonnegative$ and~$i \in \until{N}$, then the system
    asymptotically achieves multi-agent average consensus.  Moreover,
    there exists a positive MIET for each agent~$i$ given by
    \begin{align}\label{eq:miet}
      \tau_i^\text{min} = \sqrt{\frac{1}{|\NN_i|}} \left[
        \operatorname{atan}(2 \sqrt{|\NN_i|}) -
        \operatorname{atan}(\sqrt{|\NN_i|}) \right].
    \end{align}
  \end{theorem}
}

\changes{Note that the trigger in Theorem~\ref{th:dynamic2} is
  slightly different from the rules presented above in that, in
  addition to the local error~$e_i$ being reset to $0$ at each event
  triggered by agent~$i$, the internal dynamic variable~$\chi_i$ is
  reset to 1 at these times as well. The existence of a positive
  MIET~\eqref{eq:miet} makes the solution presented in
  Theorem~\ref{th:dynamic2} truly implementable on physical platforms.
  Although the solutions presented in Theorems~\ref{th:sampled-data1}
  and~\ref{th:sampled-data2} also provide a trivial
  MIET~$\tau^\text{min} = h$ guarantee based on the sampling
  period~$h$, they require perfectly synchronized executions among the
  network agents. }

Table~\ref{ta:newdetails} summarizes all proposed solutions to
Problem~\ref{pr:event-control} discussed up to this point in the
article while emphasizing their limitations. Table~\ref{ta:taxonomy} 
recalls the relevant terms.

\renewcommand{\arraystretch}{1.3}
\begin{table}[htb!]
  {
  \begin{center}
    \begin{tabular}{|ll|}
      \hline
      $x_i$ & state of agent~$i$ 
      \\
      $u_i$ & control input of agent~$i$ 
      \\
      $\hat{x}_i$ & last broadcast state of agent~$i$ 
      \\
      $f_i(\cdot)$ & triggering functions
      \\
      $L$ & Laplacian matrix 
      \\
      $z_i = (Lx)_i$ & $\sum_{j \in \NN_i} (x_i - x_j)$ 
      \\
      $\hat{z}_i = (L\hat{x})_i$ & $\sum_{j \in \NN_i} (\hat{x}_i -
      \hat{x}_j)$ 
      \\
      $\hat{\phi}_i$ & $\sum_{j \in \NN_i} (\hat{x}_i - \hat{x}_j)^2$ 
      \\
      $x_{\NN_i} = (x_i,\{ x_j \}_{j \in \NN_i})$ & state of
 agent~$i$ and neighbors
      \\
      $\hat{x}_{\NN_i} = (\hat{x}_i, \{ \hat{x}_j \}_{j \in \NN_i})
$ & last broadcast state of 
      \\ 
      & \; agent~$i$ and neighbors
      \\
      \hline
    \end{tabular}
  \end{center}
  \caption{List of terms related to Problem~\ref{pr:event-control}
    and its solutions.}\label{ta:taxonomy} 
}
\end{table}

\subsection{\textbf{Topology}}\label{se:topology}

{We have focused all our solutions so far on solving the same
  Problem~\ref{pr:event-control}.  Our discussion describes the
  different types of capabilities on the agents assumed by the
  different solutions along with their benefits and drawbacks. A
  commonality between all of them is the requirement of undirected
  communication topologies and single-integrator dynamics. From here
  on, we discuss what happens in the case of more complicated
  topologies and dynamics.}

Beyond the scenarios with undirected communication graphs considered
so far, here we extend the ideas of the article to the case where
communication topologies are directed.  The earliest works we are
aware of to address this problem are presented
in~\citep{ZL-ZC:11,ZL-ZC-ZY:12,XC-FH-AR:14,DX-JX:14}, where the
authors present either centralized or event-triggered \emph{control}
solutions only. In other words, similar to
Theorem~\ref{th:event-control}, the algorithms assume that agents have
continuous access to neighboring state information at all times. % As
% mentioned above, none of these solutions are applicable in wireless
% settings where this information must be communicated to one another.
Here, we are instead interested in event-triggered coordination
strategies similar to Theorems~\ref{th:broadcast1}
and~\ref{th:broadcast2} that not only determine when control signals
should be updated but also when communication should occur.

Here we consider communication topologies that are described by
weight-balanced directed graphs. More specifically, we say that
agent~$i$ can only send messages to its in-neighbors
$j \in \NN_i^\text{in}$ and it can only receive messages from its
out-neighbors $j \in \NN_i^\text{out}$, where the neighboring sets are
not necessarily the same. Then, consider the same type of control law
as before given in~\eqref{eq:decentralcontrol}
and~\eqref{eq:decentralcontrol2}, except each agent can now only use
information about their out-neighboring states,
\begin{align}\label{eq:directedcontrol}
  u_i(t) = - \sum_{j \in \NN_i^\text{out}} w_{ij} (\hat{x}_i - \hat{x}_j).
\end{align}

Conveniently, the closed-loop system dynamics is still given
by~\eqref{eq:decentraldynamics2} where the only difference now is that
the Laplacian~$L$ is no longer symmetric,
\begin{align}\label{eq:directeddynamics}
  \dot{x}(t) = -L \hat{x}(t) = -L(x(t) + e(t)) . 
\end{align}
The problem can now be formalized as follows.

\begin{problem}\longthmtitle{Decentralized event-triggered
    coordination on directed graphs}\label{pr:directed}
  Given a weight-balanced communication graph~$\commgraph$ and the
  closed-loop dynamics~\eqref{eq:directeddynamics}, find an
  event-trigger~$f_i(\cdot)$ for each agent~$i$ that is locally
  computable such that the sequences of times $\{ t_\ell^i \}_{\ell
    \in \integernonnegative}$ ensures multi-agent average
  consensus~\eqref{eq:averageconsensus} is achieved.
\end{problem}

Note that Problem~\ref{pr:directed} is identical to
Problem~\ref{pr:event-control} except we now consider a directed
(balanced) graph rather than an undirected one.  Just as we did in
solving Problems~\ref{pr:central} and~\ref{pr:event-control}, let us
again consider the Lyapunov function
\begin{align*}
  V(x) = \frac{1}{2}x^T L x.
\end{align*}
Given the closed-loop dynamics~\eqref{eq:directeddynamics}, we have
\begin{align*}
  \dot{V} = x^T L \dot{x} = -x^T L L (x+e).
\end{align*}
Unfortunately~$L$ no longer being symmetric causes a serious problem
because we cannot expand out this equation in a way that does not
include the in-neighbors of a given agent~$i$. As a result, using this
Lyapunov function, we are not able to find a local triggering
function~$f_i$ for agent~$i$ that only depends on the information
actually available to it (its own state and the last broadcast state
of its out-neighbors.)

Instead, let us consider the other Lyapunov
function~\eqref{eq:otherlyap}
\begin{align}\label{eq:otherlyap2}
  V(x) = \frac{1}{2} (x - \bar{x} \mathbf{1})^T (x - \bar{x}
  \mathbf{1}) ,
\end{align}
where $\bar{x} = \frac{1}{N} \sum_{i=1}^N x_i(0)$ is the average of
all initial conditions.  Then, given the closed-loop
dynamics~\eqref{eq:directeddynamics}, we have
\begin{align*}
  \dot{V} = x^T \dot{x} - \bar{x} \mathbf{1}^T \dot{x} = - x^T L
  \hat{x} - \bar{x} \mathbf{1}^T L \hat{x} = - x^T L \hat{x},
\end{align*}
where we have used the fact the  graph is weight-balanced in the last
equality.

Remarkably, a similar analysis to that used in the proof of
Lemma~\ref{le:bound} holds (after replacing~$\NN_i$
with~$\NN_i^\text{out}$ and explicitly considering weights) to yield
the following bound.

\begin{lemma}[\citep{CN-JC:16-auto}]\label{lemma:bound}
  Given~$V(x) =
  \frac{1}{2}(x-\bar{x}\mathbf{1})^T(x-\bar{x}\mathbf{1})$
  and the closed-loop dynamics~\eqref{eq:directeddynamics},
  \begin{align*}
    \dot{V} \leq \sum_{i=1}^N e_i^2 |\NN_i^\text{out}| - \frac{1}{4}
    \sum_{j \in \NN_i^\text{out}}  w_{ij} (\hat{x}_i-\hat{x}_j)^2 .  
  \end{align*}
\end{lemma}

Leveraging Lemma~\ref{lemma:bound}, it is easy to see that if we can
enforce the error of all agents to satisfy
\begin{align}\label{eq:enforce3}
  e_i^2 \leq \sigma_i \frac{1}{4 |\NN_i^\text{out}|} \sum_{j \in
  \NN_i^\text{out}} w_{ij} (\hat{x}_i - \hat{x}_j)^2 ,
\end{align}
with $\sigma_i \in (0,1)$ for all times, we have
\begin{align}\label{eq:actual}
  \dot{V} \leq \sum_{i=1}^N \frac{\sigma_i - 1}{4} \sum_{j \in
  \NN_i^\text{out}} w_{ij} (\hat{x}_i - \hat{x}_j)^2,
\end{align}
which is strictly negative for all $L\hat{x} \neq 0$. 
Letting~$\hat{x}_{\NN_i^\text{out}} = (\hat{x}_i, \hat{x}_{\NN_i^\text{out}}))$, the following
decentralized event-trigger ensures this is satisfied at all times.

\begin{theorem}\longthmtitle{Decentralized event-triggered coordination
    on directed graphs~\citep{CN-JC:16-auto}}\label{th:state-directed}
  Given {a weight-balanced communication graph~$\commgraph$
    and} the closed-loop dynamics~\eqref{eq:directeddynamics}, if the
  event times of each agent~$i$ are determined as the times when
  \begin{align}\label{eq:trigger3}
    f_i(e_i,& \hat{x}_{\NN_i^\text{out}}) \triangleq \\ 
    & e_i^2 - \sigma_i \frac{1}{4 |\NN_i^\text{out}|}
    \sum_{j \in \NN_i^\text{out}} w_{ij} (\hat{x}_i - \hat{x}_j)^2 \geq 0, \notag
  \end{align}
  then all non-Zeno trajectories of the system achieve multi-agent
  average consensus.
\end{theorem}

\begin{remark}\longthmtitle{Weight-balanced assumption} {\rm For
    implementations where the weights of the directed graph are design
    parameters, one can think of choosing them in a way that makes the
    given directed interaction topology weight-balanced. For cases
    where such choices can be made before the event-triggered
    consensus algorithm is implemented, the
    works~\citep{BG-JC:09c,AR-TC-CNH:14} present provably correct
    distributed strategies that, given a directed communication
    topology, allow a network of agents to find such weight edge
    assignments. } \oprocend
\end{remark}

In order to guarantee that the agents can converge exactly to the
average of the initial agent states, Theorem~\ref{th:state-directed}
relies on~$L$ being weight-balanced. Consequently, it is unknown if
and where the system will converge for directed graphs in general.  If
the directed graph contains a rooted spanning tree, agreement can be
reached but it is not clear to where. In other cases it may not even
be possible to reach agreement.

We have now discussed different types of static or fixed communication topologies,
but one could easily imagine scenarios with both
time-varying~\citep{WZ-HL-ZJ:16} or state-dependent interaction graphs
for which some modified algorithms may need to be developed. In the case
of state-dependent interaction graphs, an additional challenge that
must be addressed is how to preserve connectivity of the network while
performing the primary consensus task~\citep{HY-PJA:12,YF-GH:15,XY-JW-DVD-KHJ:17,PY-LD-ZL-ZG:16}.

\subsection{\textbf{Dynamics}}\label{se:dynamics}

So far, we have only considered the simple single-integrator
dynamics~\eqref{eq:dynamics}. While these simple dynamics are
certainly useful for virtual states (e.g., a temperature estimate) or
in demonstrating the ideas of event-triggered consensus in general,
this might be too limited in cases where states correspond to physical
quantities. % Consequently, the consideration of more
% complicated dynamics to make them applicable to more realistic
% scenarios are natural extensions.
In this section we begin by discussing double-integrator dynamics
before moving onto general linear dynamics. We note here that as we
generalize the dynamics we also consider goals beyond average
consensus. More specifically, depending on the dynamics, static
average consensus may not be possible, in which case we will instead aim for
synchronization.

\subsubsection{Double-integrator systems} \label{sec:DI}
%
%{\color{red}
%The latter is concerned with generalizations of event-triggered consensus algorithms for consensus of second-order systems, systems with general linear dynamics (refer to Sections \ref{sec:DI} and \ref{sec:Linear}, respectively), and consensus on different types of information such as plan assignments \cite{MA-JPH:08}, \cite{HLC-LB-JPH:09}, \cite{SSP-JR-HLC-JPH-MV-JV:10}.
%}

Let us consider the case where the state of agent~$i$ is denoted
by~$x_i = (r_i,v_i) \in \real \times \real$ with double-integrator
dynamics,
\begin{align}\label{eq:dynamicsdouble}
  \dot{r}_i(t) &= v_i(t) , \\
  \dot{v}_i(t) &= u_i(t) . \notag
\end{align}
Then, it is known that the distributed controller
\begin{align}\label{eq:double-ideal}
  u_i^*(x) = - \sum_{j \in \NN_i} (r_i - r_j + \gamma(v_i - v_j) )
\end{align}
with $\gamma > 0$ drives the states of all agents to a consensus
trajectory~\citep{WR-EA:07,WR:08b}.  Different from consensus with
single-integrator dynamics, where the agents are to reach a steady
state, the consensus problem of double-integrator systems is rather to
synchronize the outputs.  More specifically, we say that a system
asymptotically achieves synchronization if
\begin{align}\label{eq:synchronization}
  \lim_{t \rightarrow \infty} (x_i(t) - x_j(t)) = 0. 
\end{align}
Note in particular that different from the definition of consensus we
used for single-integrator dynamics~\eqref{eq:averageconsensus}, we
only require $x_i(t)$ and $x_j(t)$ to tend together as time goes to
infinity rather than both going to the same stationary point.
This is formalized in Theorem~\ref{th:DI}.

\begin{theorem}\longthmtitle{Continuous controller
    (double-integrators)~\citep{WR:08b}}\label{th:DI}
  Given {a connected, undirected graph~$\commgraph$ and} the
  dynamics~\eqref{eq:dynamicsdouble}, if all agents implement the
  control law~\eqref{eq:double-ideal}, then the system asymptotically
  achieves synchronization~\eqref{eq:synchronization}.
\end{theorem}

Implementing~\eqref{eq:double-ideal} requires agents to have
continuous information about one another. Thus, we are interested in
an event-triggered implementation of the ideal control
law~\eqref{eq:double-ideal} to relax this
requirement. Letting~$\{t_\ell^i\}_{\ell \in \integernonnegative}$ be
the sequence of event-times for agent~$i$, we let
\begin{align}\label{eq:xhatDI}
  \hat{x}_i(t) = x_i(t_\ell^i) = (r_i(t_\ell^i),v_i(t_\ell^i)) \text{
  for } t \in [t_\ell^i,t_{\ell+1}^i) 
\end{align}
be the last broadcast state of agent~$i$. Then, at any given time~$t$,
agent~$i$ only has access to the last broadcast state~$\hat{x}_j(t)$
of its neighbors~$j \in \NN_i$ rather than exact states~$x_j(t)$.

The distributed event-triggered controller we consider is then given by
\begin{align}\label{eq:DIcontrol}
  u_i(t) &= 
  \\
  -& \sum_{j \in \NN_i} (\hat{r_i} + (t-t_\ell^i) \hat{v}_i - \hat{r}_j
     - (t-t_\ell^i) \hat{v}_j + \gamma(\hat{v}_i - \hat{v}_j)) \notag, 
\end{align}
for $t \in [t_\ell^i, t_{\ell+1}^i)$. It should be noted that this
controller utilizes a first-order-hold (FOH) instead of a ZOH for
the~$r_j$ components of the state. We now define two different errors
\begin{align*}
  e_{r,i}(t) &= \hat{r}_i(t) + (t-t_\ell^i)\hat{v}_i(t) - r_i(t), 
  \\
  e_{v,i}(t) &= \hat{v}_i(t) - v_i(t).
\end{align*}
Then, defining the stack vector of the
error~$e = \left[ \begin{array}{cc} e_r^T & \gamma e_v^T \end{array}
\right]^T$,
the closed-loop dynamics of the controller~\eqref{eq:DIcontrol} is
given by
\begin{align}\label{eq:DI-closed}
  \left[ \begin{array}{c} \dot{r}(t) \\ \dot{v}(t) \end{array} \right] %= \Gamma \left[ \begin{array}{c} \hat{r}(t) \\ \hat{v}(t) \end{array} \right] 
  = \Gamma \left[ \begin{array}{c} r(t) \\ v(t) \end{array} \right] - \left[ \begin{array}{cc} 0 & 0 \\ L & L \end{array} \right] e(t),
\end{align}
where
\begin{align}
  \left.
	\begin{array}{l l}
          \Gamma=\begin{bmatrix}
            \mathbf{0}_{N \times N} & I_N \\
            - {L} & -\gamma {L}
\end{bmatrix} 
   \end{array}  \right. . \label{eq:Gamma}
\end{align}
The problem can now be formalized as follows.

\begin{problem}\longthmtitle{Decentralized event-triggered
    coordination with double-integrator dynamics}\label{pr:DI} 
  Given a connected, undirected graph~$\commgraph$ and the closed-loop
  dynamics~\eqref{eq:DI-closed}, find an event-trigger~$f_i(\cdot)$
  for each agent that is locally computable such that the sequences of
  times~$\{t_\ell^i\}_{\ell \in \integernonnegative}$ ensures
  synchronization~\eqref{eq:synchronization} is achieved.
\end{problem}

This problem is first addressed
in~\citep{DX-SX-ZL-YZ:15,MC-FX-LW:14} where state-dependent
triggering rules are developed; however, these algorithms rely on continuous
state information about neighbors. The first work we are aware of that solves
Problem~\ref{pr:DI} where information is only shared at event times
is~\citep{GSS-DVD-KHJ:13}, which proposes the time-dependent triggering
threshold presented next.

\begin{theorem}\longthmtitle{Decentralized event-triggered
    coordination with double-integrator
    dynamics~\citep{GSS-DVD-KHJ:13}}\label{th:DI-time} 
  Given {a connected, undirected graph~$\commgraph$ and} the
  closed-loop dynamics~\eqref{eq:DI-closed}, if the event times of
  each agent~$i$ are determined as the times when
  \begin{align}\label{eq:DI-trigger-time}
    f_i(e_i,t) \triangleq \Big\Vert \left[ 
      \begin{array}{c}
        e_{r,i}
        \\
        \gamma e_{v,i}
      \end{array}
    \right] \Big\Vert - (c_0 + c_1 e^{-\alpha t}) = 0,
  \end{align}
  with constants~$c_0,c_1 \geq 0$ and $c_0 + c_1 > 0$, then all
  non-Zeno trajectories of the system reaches a neighborhood of
  consensus upper-bounded by
  \begin{align*}
    r = \TwoNorm{L} \sqrt{2N} c_0 c_V / \lambda_3(\Gamma) ,
  \end{align*}
  where~$c_V > 0$ is related to the graph~$L$. Moreover, if $c_0 > 0$
  or $0 < \alpha < \lambda_3(\Gamma)$, then the closed-loop system
  does not exhibit Zeno behavior.
\end{theorem}
%\begin{pf}
%  In the appendix. \hfill \qed
%\end{pf}

Similar to the result of Theorem~\ref{th:time}, this solution uses a
triggering function whose threshold depends on time rather than
state. However, here we have also considered the use of a first-order
holder controller between events. More specifically, this
problem-solution pair has double-integrator dynamics, an undirected
interaction graph, events that trigger broadcasts and control updates,
triggers that are evaluated continuously, and trigger thresholds that
are time-dependent.

Also similar to the algorithm presented in Theorem~\ref{th:time}, this
solution has design parameters~$c_0,c_1,$ and~$\alpha$ that can be
tuned to balance performance and efficiency. However, there are also
limits to how these parameters can be tuned. In particular, we recall
that there are only two ways to guarantee Zeno behavior. The first is
to set $c_0 > 0$, but this sacrifices being able to guarantee
convergence all the way to the exact consensus state. The second is to
set~$\alpha < \lambda_3(\Gamma)$, but this requires global knowledge
of the entire network structure to be able to compute.

We have only considered undirected topologies here there are indeed
works that have considered directed topologies as well.  For directed
topologies, even when the graph contains a spanning tree, the ideal
controller~\eqref{eq:double-ideal} (with continuous communication and
actuation) is only guaranteed to drive the system to a consensus state
if~$\gamma$ is sufficiently large. We omit the details and instead
refer the interested reader to~\citep{GSS-DVD-KHJ:13}.

Other works addressing event-triggered consensus of double-integrator
systems include \citep{HY-YS-HZ-HS:14,DX-SH:13,XY-DY:13,NM-XL-TH:15}
with various modifications. The work by \cite{DX-SH:13} consider the
case of heterogeneous communication networks. In this case the
positions and velocities are shared using different communication
graphs. % $\mathcal{G}_r$ and $\mathcal{G}_v$, respectively
The work~\cite{NM-XL-TH:15} discusses the leader-follower consensus
problem with double-integrator systems and \citep{EG-YC-DWC:16}
addresses decentralized event-triggered consensus of second-order
systems in the presence of communication imperfections. In particular,
this reference considered the presence of communication delays and
packet dropouts using a broadcasting style of
communication. Discrete-time systems have also been explored in this
context and we refer the interested reader
to~\citep{XC-FH:12,WZ-HP-DW-HL:17,XY-DY:13,XY-DY-SH:13}.

\subsubsection{Linear systems}   \label{sec:Linear}

The event-triggered approach has also been extended to consider more
general dynamics than double-integrator models.  Here we discuss in
detail the synchronization problem for a homogeneous group of~$N$
agents or subsystems with linear dynamics.

Letting~$x_i \in \real^n$ denote the state of agent~$i$, we consider
homogeneous linear dynamics
\begin{align}\label{eq:linear}
  \dot{x}_i(t) = Ax_i(t) + Bu_i(t) ,
\end{align}
where the pair~$(A,B)$ is controllable.  The objective of the
consensus or synchronization problem of~\eqref{eq:linear} is to drive
the state of each system~$x_i$ to a common response or trajectory,
that is, the corresponding elements of each agent's state need to
converge to a single trajectory.

The synchronization of multi-agent systems with linear dynamics and
assuming continuous communication has been extensively studied,
e.g.,~\cite{ZL-ZD-GC-LH:10,ZL-ZD-GC:11,CQM-JFZ:10,WR:08,LS-RS:09,YS-JH:12,SET:08,SET:09,JHS-HS-JB:09}.
It is known~\citep{CQM-JFZ:10,EG-YC-DWC:14} that the distributed
controller
\begin{align}\label{eq:LDinputs}
  u_i(t)= c F \sum_{j\in\mathcal{N}_i}(x_i(t)-x_j(t)) 
\end{align}
with~$c>0$ ensures that the system achieves synchronization under some
suitable conditions on the matrix~$F$.
%%
%\marginJC{What's the dimension of $F$?}
%\margin{dimension of $u_i$ times $n$}
%%
This is formalized next.

\begin{lemma}[\citep{EG-YC-DWC:14}]\label{th:NS_LD}
  Given the dynamics~\eqref{eq:linear}, if all agents implement the
  control law~\eqref{eq:LDinputs} and the matrices
  $(A + c \lambda_j(L) B F)$ are Hurwitz for $j = 2,3, \dots, N$, then
  the system asymptotically achieves synchronization.
\end{lemma}

Unfortunately, the condition of Lemma~\ref{th:NS_LD} requires checking
several matrices that are functions of every nonzero eigenvalue of the
Laplacian~$L$ of the communication graph. This is a  difficult
condition to check in general as it requires knowledge of all
eigenvalues of~$L$ in order to compute the eigenvalues of the matrices
$(A + c \lambda_j(L) B F)$, for $j = 2,3, \dots, N$. Less restrictive
conditions involve the design of the consensus protocol parameter
using only the smallest non-zero eigenvalue of the Laplacian matrix,
$\lambda_2$, or an \emph{estimate} of such eigenvalue. The following
result gives one way of designing the controller~\eqref{eq:LDinputs}
to satisfy the conditions of Lemma~\ref{th:NS_LD}.

\begin{theorem}\longthmtitle{Continuous controller (linear
    dynamics)}\label{th:ConsLin}
  Given {a connected, undirected graph~$\commgraph$ and}
  ($A,B$) controllable, if all agents implement the protocol
  \eqref{eq:LDinputs}, with
  \begin{align}
    F&=-B^TP, \label{eq:F}\\
    c&\geq1/\lambda_2, \label{eq:c}
  \end{align}
  where~$P$ is the unique solution to
  \begin{align}
    PA+A^TP-2PBB^TP+2\alpha P<0, \label{eq:LMI}
  \end{align}
  then the system ensures synchronization~\eqref{eq:synchronization} is achieved.
%  Furthermore, the Lyapunov function defined by
%  $V=x^{T}\hat{\mathcal{L}}x$ has a time derivative along the
%  trajectories of~\eqref{eq:linear} with inputs~\eqref{eq:LDinputs}
%  given by $\dot{V}=x^T\bar{\mathcal{L}}x$.
\end{theorem}
% Implementing~\eqref{eq:LDinputs} in a digital setting is not possible
% since it requires all agents to have continuous access to the state of
% their neighbors. 
Next, we turn our attention to seek event-triggered implementations of
the ideal control law~\eqref{eq:LDinputs}. One of the earliest works
to address this problem is~\citep{WZ-ZJ-GF:14}, which
considers a digital implementation of the ideal
controller~\eqref{eq:LDinputs}
\begin{align}\label{eq:LDinputsZOH}
  u_i(t)=cF\sum_{j\in\mathcal{N}_i}(\hat{x}_i(t) - \hat{x}_j(t)), 
\end{align}
where
\begin{align}
  \hat{x}_i(t) = x_i(t_\ell^i) \text{ for } t \in [t_\ell^i, t_{\ell+1}^i )
\end{align}
denotes the state of agent~$i$ at its last event
time. \cite{WZ-ZJ-GF:14} then propose a simple constant threshold
event-triggering algorithm where agent~$i$ broadcasts its state to its
neighbors whenever its error~$e_i(t) = \hat{x}_i(t) - x_i(t)$ exceeds
a fixed threshold.

Unfortunately, this algorithm faces some issues resulting from the
fact that the closed-loop system may actually be unstable. More
specifically, the algorithm proposed in~\citep{WZ-ZJ-GF:14} generally
provides poor performance in terms of reducing control updates and may
lead to Zeno behavior when the agents have unstable dynamics, i.e.,
when one or more eigenvalues of the matrix $A$ have positive real
parts. Under this scenario the exponential divergence of the states of
the agents causes the error $e_i$ to grow greater than the fixed
threshold used in \cite{WZ-ZJ-GF:14} faster and faster ultimately
leading to the undesired Zeno behavior.

The event-triggered consensus of linear systems using ZOH
implementations was also addressed in \cite{GG-LD-QH:14}. In this
reference the agents measure their state and evaluate their
event-triggered conditions periodically, at every $h$ time
units. However, the convergence conditions expressed
in~\cite{GG-LD-QH:14} require explicit knowledge of the Laplacian
matrix, which is an impediment for decentralized implementation.
% hence, the design of the consensus protocol and the selection of the
% sampling period $h$ cannot be performed in a decentralized manner
% following this approach.

The early work by~\cite{TL-DJH-BL:12} provides an event-triggered
strategy to reduce communication of a class of linear systems without
explicitly addressing Zeno behavior.  The works~\cite{OD-JL:12a}
and~\cite{OD-JL:12b} study event-triggered synchronization of linear
systems using a ZOH implementation with constant thresholds. The
work~\cite{ZZ-FH-LZ-LW:14} uses a ZOH method but restricts only
actuation updates, while continuous communication is still required
for the agents to determine the triggering instants. The
work~\cite{BZ-XL-TH-GC:15} addresses leader-follower consensus problem
of linear systems but, similar to the previous reference, the event
conditions require continuous state information from neighbors and
from the leader, which limits the application of this approach for
reducing communication frequency.

Fortunately, some new frameworks and algorithms have recently been
developed to overcome this
problem~\citep{NM-XL-TH:15b,DD-ZW-BS-GW:15,DY-WR-XL:14,DY-XL-WC:15,DY-WR-XL-WC:16,EG-YC-DWC:14,QL-ZW-XH-DHZ:15,EG-YC-DWC:17,HZ-GF-HY-QC:14,CDP:13}. Recent
contributions have relied on model-based or estimation approaches in
order to address event-triggered coordination, and to
improve the performance of the multi-agent coordination system in
terms of asymptotic convergence and reduction of generated
events. Some of these approaches rely on both, sensing the states of
neighbors and transmitting the control inputs, while other approaches
only assume broadcasting capabilities.

More specifically, we now make one big change to the definition
of~$\hat{x}_i(t)$. Until now we have always treated~$\hat{x}_i(t)$ as
a piece-wise constant value that only changed when agent~$i$ triggered
an event. Instead, with a slight abuse of notation, we now consider it
as a time-varying estimate of the state of agent~$i$ with dynamics
\begin{align}\label{eq:LDagentsModel}
  \dot{\hat{x}}_i(t)=A\hat{x}_i(t)+B\hat{u}_i(t),
\end{align}
where
\begin{align}
  \hat{u}_i(t) = u_i(t_\ell^i) \text{ for } t \in [t_\ell^i,t_{\ell+1}^i)
\end{align}
is the control input used by agent~$i$ at its last event time. Then,
when an event is triggered by agent~$i$ at time~$t_\ell^i$, it sends
its control input~$u_i(t_\ell^i)$ in addition to its current
state~$x_i(t_\ell^i)$. With this information, any
neighbor~$j \in \NN_i$ can then propagate the
estimate~\eqref{eq:LDagentsModel} forward in time.

Given these new model-predictive estimates, we redefine the control
law in the same way as~\eqref{eq:LDinputsZOH} but
% \begin{align}\label{eq:LDinputsmodel2}
%   u_i(t)=cF\sum_{j\in\mathcal{N}_i}(\hat{x}_i(t) - \hat{x}_j(t)), 
% \end{align}
with~$\hat{x}_j(t)$ now given by~\eqref{eq:LDagentsModel}, rather than
being piece-wise constant.  The closed-loop system dynamics
of~\eqref{eq:linear} with these control inputs is then given by
%\begin{align}\label{eq:LD-closed}
%\dot{x}_i(t) = A x_i(t) + c B F (Lx(t))_i + c B F (Le(t))_i
%\end{align}
\begin{align}\label{eq:LD-closed}
  \dot{x}(t) = I_N \otimes A x(t) + c L \otimes B F x(t) + c L \otimes
  B F e(t) . 
\end{align}
The problem can now be formalized as follows.

\begin{problem}\longthmtitle{Decentralized event-triggered
    coordination with linear dynamics}\label{pr:linear} 
  Given the closed-loop dynamics~\eqref{eq:LD-closed}, find an
  event-trigger~$f_i(\cdot)$ for each agent that is locally computable
  such that the sequences of times~$\{t_\ell^i\}_{\ell \in
    \integernonnegative}$ ensures
  synchronization~\eqref{eq:synchronization} is achieved.
\end{problem}

% Similarly, the control input $\hat{u}_j(t)$ in
% \eqref{eq:LDagentsModel} can be updated at the event time instants
% $t_k^j$ if it is assumed that the agent also transmits the value
% $u_j(t_k^j)$ along with the sampled state $x_i(t_k^j)$. Different
% approaches exist that provide different ways to update the model
% control input $\hat{u}_j(t)$ that depend on sensing capabilities and
% the communication topology, i.e. directed or undirected graphs. In
% other papers, this control input is completely disregarded and the
% model equations take the simplified form
%\begin{align}
%	\dot{\hat{x}}_j(t)=A\hat{x}_j(t).    \label{eq:LDagentsModelSim}
%\end{align}
%This method relies in the property that the control inputs $u_i(t)$
% for $i=1,...,N$ will vanish as time progresses (under a properly
% designed event-triggered consensus algorithm). Thus, better
% estimates of the state of the real system \eqref{eq:LDagents} are
% obtained using the models \eqref{eq:LDagentsModelSim} as time
% elapses because the control input eventually disappears from
% \eqref{eq:LDagents}. In the case where simplified model form
% \eqref{eq:LDagentsModelSim} is implemented, the required information
% to be transmitted is also reduced since the agents only need to
% broadcast their states $x_j(t_{k_j})$ but not their control inputs
% $u_j(t_{k_j})$.

Leveraging the result of Theorem~\ref{th:ConsLin}, a state-dependent
triggering rule to solve Problem~\ref{pr:linear} is proposed
in~\citep{EG-YC-DWC:14} based on the Lyapunov
function~$V = x^T \LL x$, where $\LL = L \otimes P$ and~$P$ is defined
by~\eqref{eq:LMI}.  This result is formalized next.

\begin{theorem}\longthmtitle{Decentralized (state-dependent)
    event-triggered 
    coordination with linear dynamics~\citep{EG-YC-DWC:14}}\label{th:Decent}
  Given {a connected, undirected graph~$\commgraph$ and} the
  closed-loop dynamics~\eqref{eq:LD-closed}, if the event times of
  each agent~$i$ are determined as the times when
  \begin{align}\label{eq:DecCond}
    f_i(e_i,\hat{x}_{\NN_i}) \triangleq \delta_i -
    \sigma_i \hat{z}_i^T \Theta_i \hat{z}_i \geq 0,
  \end{align}
  where
  \begin{align*}
    \hat{z}_i&=\sum_{j\in\mathcal{N}_i}\left(\hat{x}_i-\hat{x}_j\right),
    \\
    \Theta_i& =(2c_2-b_i|\NN_i|(c_2-c_1))PBB^{T}P,
    \\
    \delta_i&=2(c_2-c_1)|\NN_i|\hat{z}_i^TPBB^TPe_i + |\NN_i|
    e^T_iPBB^TPe_i
    \\
    &\times \left[2c_1|\NN_i|(1+b_i)+\frac{c_2-c_1}{b_i}\right.
    +\left.c_1(N-1)\left(b_i+\frac{3}{b_i}\right)\right],
  \end{align*}
  with constants~$c_1\geq 1/\lambda_2$, $c_2>0$, and
  $0<b_i<\frac{2c_2}{|\mathcal{N}_i|(c_2-c_1)}$ for $c_2>c_1$, or $b_i > 0$
  otherwise, then all non-Zeno trajectories of the system
  asymptotically achieve synchronization.
\end{theorem}

Theorem~\ref{th:Decent} guarantees asymptotic synchronization of the
agents with linear dynamics along all its non-Zeno
trajectories. Unfortunately, this result does not guarantee the
exclusion of Zeno behavior. Following our discussion in
Section~\ref{se:zeno}, there are now several methods that can be used
to address this issue. For instance, in~\citep{EG-YC-DWC:14} a small
fixed parameter is included in the trigger function \eqref{eq:DecCond}
to avoid Zeno behaviors.  While the modified algorithm is able to
ensure Zeno behavior do not occur, the price to pay is that it can no
longer guarantee convergence exactly to a synchronized state, but
rather to a neighborhood around it. We omit the details but refer the
interested reader to~\citep{EG-YC-DWC:14}.

For completeness, the next result from~\citep{EG-YC-DWC:17b} also
solves Problem~\ref{pr:linear} using a time-dependent triggering
function rather than a state-dependent one as in
Theorem~\ref{th:Decent}.

\begin{theorem}\longthmtitle{Decentralized (time-dependent) event-triggered
    coordination with linear dynamics~\citep{EG-YC-DWC:17b}}\label{th:LD-time}
  Given {a connected, undirected graph~$\commgraph$ and} the
  closed-loop dynamics~\eqref{eq:LD-closed}, if the event times of
  each agent~$i$ are determined as the times when
  \begin{align*}
    f_i(e_i,t) \triangleq \TwoNorm{e_i} - c_1 e^{-\alpha t} = 0,
  \end{align*}
  with constants $c_1, \alpha > 0$, then all non-Zeno trajectories of
  the system asymptotically achieve synchronization. Moreover, there
  exists $\lambda^*$ such that for $\alpha < \lambda^*$, the
  closed-loop system does not exhibit Zeno behavior.
%The parameter $\hat{\lambda}$ is such that $ \left\|\textbf{e}^{\hat{A}t}\right\|\leq\hat{\beta}\textbf{e}^{-\hat{\lambda} t}$, for $\hat{\beta}>0$ where $\hat{A}=I_{N-1}\otimes A +cJ_{2:N}\otimes BF$.
%Furthermore, the agents do not exhibit Zeno behavior and the inter-event times $t_{k_i+1}-t_{k_i}$ for every agent $i=1,...,N$ are bounded by the \textit{positive} time $\tau$, that is
%\begin{align}
%  \tau\leq t_{k_i+1}-t_{k_i}  \label{eq:tk}
%\end{align}	
%where 
%\begin{align}
%\tau=\frac{\ln(1+\beta/K_\tau)}{\left\|A\right\|+\hat{\lambda}}  \label{eq:tausol}
%\end{align}
%and the parameter $K_\tau>0$ is a function of the design and the connectivity parameters.
\end{theorem}

Similar to the result of Theorem~\ref{th:time}, Zeno behavior can be
avoided if~$\alpha$ is chosen small enough, where the critical
value~$\lambda^*$ again depends on global information. The interested
reader is referred to~\citep{EG-YC-DWC:17b}. Similar problems in the
presence of input saturation are considered
in~\citep{XL-CD-PL-DY:16,BZ-XL-TH-HL-GC:16}.

We stop our discussion here at homogeneous agents with linear
dynamics. Indeed there are also some relevant works that address the
case of heterogeneous agents with linear
dynamics~\citep{FZ-ZH-YY-JW-LL-JP:17} as well. Additionally, recent
works have also considered nonlinear dynamics but we omit the details
and instead refer the interested reader to the
works~\citep{YX-YD:13,XZ-LL-GF:15,HL-GC-LX:16,HL-GC-TH-WZ-LX:16,AH-JC:17,DL-DVD-MB-KHJ:16}. Instead,
some works consider nonlinear control inputs with simpler dynamics to
achieve finite-time or fast consensus~\citep{MG-DVD:13,XZ-MC-LW:15}.

% {\color{red}~\cite{XZ-PS-CL-CY-WG:15} - linear dynamics, undirected
% broadcast, considers a cost function, considers delays}

\subsection{Uncertainty}\label{se:uncertainties}

Throughout this article we have assumed there are no disturbances or
uncertainties of any kind, which clearly idealizes many instances of
the problems we are interested in. Consequently, it is clearly
important to determine how robust the algorithms we have discussed so
far are in the presence of different sources of uncertainty, and how
the algorithms might need to be modified to accommodate them.  While
we do not go into the same level of detail for these algorithms as we
have done in the rest of the article, we provide a brief discussion of
the technical issues and challenges on this front.

\changes{As discussed in Section~\ref{se:zeno}, guaranteeing the
  existence of a positive MIET is crucial in ensuring that a proposed
  solution can actually be implemented.  Interestingly, even if an
  event-triggered controller guarantees a positive MIET when
  disturbances do not exist, it is possible that arbitrarily small
  disturbances are enough to void this guarantee. The robustness of
  the existence of a positive MIET is therefore a critical
  property. We refer to~\cite{DPB-WPMHH:14} for a detailed discussion
  on the notion of robustness of event-triggered controllers in the
  presence of disturbances and, in particular, the notion of
  \emph{local event-separation} (guaranteeing non-Zeno executions when
  disturbances are not present) and the stronger notions of
  \emph{semi-global} and \emph{global event-separation} (with
  event-triggered controllers being robust against disturbances).  The
  implementation details when designing event-triggered communication
  and control algorithms highlight the importance of this topic. While
  this article has stopped just short of this technical discussion,
  future works on these topics should be mindful of these
  implementation details to ensure the solutions can be practically
  implemented.  }

\subsubsection{Quantization}
\changes{Beyond state disturbances, imperfect wireless communication
  mechanisms present an additional set of challenges. } As a majority
of this article assumes messages are shared wirelessly, quantization
of transmitted information is a natural issue that must be
addressed. This problem was first studied by \cite{EG-YC-HY-PA-DC:13}
using uniform quantizers. In general, an algorithm which would
guarantee asymptotic convergence to the global initial average when
non-quantized information is transmitted by each agent, is instead
only able to converge to a bounded region around the initial average
in the case where uniform quantizers are implemented at each
node. This type of result is also known as \emph{practical
  consensus}. The states of the agents satisfy $\lim_{t\rightarrow
  \infty} ||x_i(t)-\bar{x}||\leq c$, where $c$ is a constant which
depends on the size of the quantization step and
$\bar{x}=\frac{1}{N}\sum_{i=0}^Nx_i(0)$.

The more recent work~\citep{ZZ-LZ-FH-LW:15} considers both uniform and
logarithm quantizers. In the case of logarithmic quantizers,
asymptotic convergence to the initial average is still not guaranteed
in general. The difference between any two states is still bounded,
but now this bound also depends on the value of the initial
average. For instance, if the initial average happens to be equal to
zero, then asymptotic convergence to the initial average is
achieved. Compared to uniform quantizers, the use of logarithmic
quantizers has been of significant advantage in stabilization problems
since the quantization error diminishes as the signal to be quantized
tends to zero. However, in consensus problems of single-integrator
systems, the steady state value of the overall system is in general
not equal to zero. Hence, it is expected that the results shown by
\cite{ZZ-LZ-FH-LW:15} do not guarantee asymptotic convergence and the
bounds on the state disagreement depend on the value of the initial
average. Some more recent works that study further extensions such as
dynamic quantization or self-triggered mechanisms are proposed
in~\citep{CDP-PF:13,HL-GC-TH-ZD-WZ-LG:16,XY-JW-KHJ:16,DS-PT-CDP:18}.

\subsubsection{Communication delays and packet dropouts}
In addition to quantization, specifically considering wirelessly
networked systems introduces many new sources of uncertainty in the
form of delays and packet drops. Consequently, some of these issues
have already been studied when the algorithms were first
developed. Early in the development of event-triggered consensus
algorithms,~\citep{GSS-DVD-KHJ:13} analyzes the presence of
communication delays in the consensus of single-integrator systems. It
was shown in that reference that the closed-loop overall system using
event-triggered controllers is Input-to-State Stable (ISS) with
respect to the state errors introduced by the event-triggered
controllers and practical consensus was demonstrated in the presence
of delays bounded by a function of the largest eigenvalue of the graph
Laplacian. Consensus of discrete-time single-integrator systems with
communication delays was studied in~\citep{LL-DH-SX:14}. The
work~\citep{EG-YC-DWC:16} provides an approach for consensus of
double-integrator systems using a time-dependent threshold for systems
with non-consistent packet dropouts and delays. Non-consistent packet
dropouts means that a packet of information broadcasted by a given
agent $i$ may be received by all, some, or none of the intended
recipients $j$ such that $i\in \mathcal{N}_j$. Similarly, the
communication delay associated to a given broadcast message can be
different in general to every receiving agent, given that the message
is successfully received. {This is studied in~\cite{VSD-MA-WPMHH:17}
for general linear systems where a dynamic event-triggered
coordination strategy is developed that guarantees average consensus
but requires global information to determine whether the system will
converge or not.}

\section{Applications of Event-Triggered
  Consensus}\label{se:applications}

{Here we provide examples of both direct and indirect
  applications of the various algorithms discussed in the article.
  Our presentation is not meant to be exhaustive of every area, but
  rather serve as an initial point of reference for interested readers
  that seek to employ the results presented above. Our focus is not on
  multi-agent consensus per se, but in works that have addressed the
  need for event-triggered coordination in problems that involve
  consensus.}

\paragraph*{Formation control.} The connections between consensus and
some practical problems such as formation control of groups of
vehicles have already been long established before event-triggered
ideas became popular~\citep{WR-RWB:08,WR-EA:07}. {This problem is relatively straightforward
in the centralized case where all agents know the desired shape and location of the
final formation. However, in the decentralized version of the problem, each vehicle 
may know the desired formation shape but the location of the formation needs to be
negotiated and agreed upon by the distributed agents~\citep{WR-RWB-EMA:07}. 

Letting~$p_i \in \real^n$ represent the position of agent~$i$ in some
space, the goal is to drive all agents~$i \in \until{N}$ such
that~$p_i(t) \rightarrow \bar{p} + b_i$, where~$\bar{p}$ is the
average position of all agents and~$b_i$ represents the desired
relative displacement of agent~$i$ defining its position in the
formation with respect to the center. Since the vectors~$b_i$ are
constant, the agents then simply need to perform average consensus on
the virtual state~$x_i = p_i - b_i$ to agree upon the center of the
formation.}  Early works generally assumed that agents have
continuous, or at least periodic, access to information about their
neighbors. Consequently, several groups of researchers have considered
applying the event-triggered coordination ideas to consensus-based
{formation control} algorithms to relax this
requirement~\citep{XC-ZP-GW-AR:17,CN-GJP:16,AA-DL-DVD-KHJ:17}.

\paragraph*{Leader-tracking.} Similar to formation control, another
popular application is to actually try to directly control the entire
group of agents using a small subset (or even just a singular) of
agents to drive the rest of the group. In the {leader-tracking (or
  leader-follower)} problem there exists (at least) one particular
agent, called the leader (or pinned agent), which acts independently
from all other agents' states. The rest of the agents are referred to
as the followers and they implement some form of consensus algorithm
such that they essentially `follow' the leader(s), perhaps while
maintaining a specified formation.  

{For simplicity, consider a single leader with identity~$i =
  0$ which aims to lead a group of agents with identities~$i =
  \until{N}$. The leader is free to move or be controlled by a
  user. Not all agents will have access to direct information about
  the leaders' motion, and hence they implement average consensus to
  propagate it throughout the network and be able to follow the
  leader.  % Since the leader is the only one
  % not actively trying to maintain the average state of all agents,
  % it
  % essentially has the power to influence the entire network.
  Depending on constraints such as the maximum speed or acceleration
  of the leader, various results can be established regarding the
  global behavior and performance of the system.  } Many groups of
researchers have looked at applying event-triggered coordination to
many different variations of this leader-tracking problem. For
instance, single- and double-integrator dynamics are considered
in~\citep{KL-ZJ-GX-RX:16}, homogeneous linear dynamics
in~\citep{CT-KZ-JMS-WED:14,YC-VU:16,WZ-ZJ:15}, heterogeneous linear
dynamics in~\cite{EG-YC-DWC:17b}, and nonlinear dynamics
in~\cite{AA-FA-DL-GS-DVD-MB-KHJ:15,XZ-LL-GF:15,HL-GC-LX:17}. In
addition to considering different types of dynamics, other groups have
considered discrete-time systems~\citep{MC-LZ-HS-CL:15}. Finally, even
more variations can be considered by imposing different types of
constraints on the problems or solutions as discussed throughout this
article. Examples include the addition of uncertainties/disturbances,
specific goals of a network of leader-followers (e.g., containment
control), dynamic topologies, and the different types of triggering
mechanisms discussed in this
article~\citep{HL-XL-TH-WZ:15,HL-GC-XL-TH:16,WX-DH-LL-JC:17,NM-XL-TH:15}.

\paragraph*{Distributed estimation.} A popular indirect application of
consensus algorithms in general is distributed state estimation, see
e.g.~\citep{LZ-ZW-DZ:17}. Consensus protocols can then be used to
allow distributed agents to communicate and agree on a common state
estimate; however, this generally assumes periodic communication among
the agents. {More specifically, distributed agents are sharing
  new samples with one another at all times to maintain both a good
  and consistent estimate of the quantity of interest. However, this
  can be wasteful in general, especially if new samples being shared
  are not providing much new information.  } Instead, applying
event-triggered coordination to these algorithms can help reduce the
amount of communication required by a network to maintain a state
estimate. In~\citep{GB-LC-DS:16}, each node of the communication
network implements a local Kalman filter and shares this information
with neighbors to achieve consensus. The implementation of
event-triggered communication strategies restricts inter-agent
communication. In this case, it is necessary to evaluate the
difference between probability density functions (PDF) at different
time instants. This is achieved by applying the Kullback-Leibler
divergence metric to the current local PDF and the last transmitted
PDF. The result of this operation is compared against a time-dependent
threshold of the form $c_0+c_1e^{-\alpha t}$. {Similar to the
  ideas of this article, the agents only share their new information
  with neighbors if the new information is different enough from what
  is currently estimated according to the time-varying
  threshold. ~\cite{MO-DI-NA-SM:18} pursue similar ideas in the
  context of cooperative localization, with agents only sending
  measurements to neighbors when the expected innovation for state
  estimation is high.} Instead, the work~\citep{QL-ZW-XH-DHZ:15b}
considers a simpler positive threshold parameter to dictate events.

\paragraph*{Clock synchronization.} Clock synchronization is a
particular problem of interest requiring agents in a network to
synchronize their imperfect clocks via communication. This problem has
been addressed in the past using consensus algorithms and assuming
continuous or periodic communication such as in~\citep{LS-FF:11}
and~\citep{RC-SZ:14}.  {Specifically, letting~$t \in
  \realnonnegative$ represent the true global time, we define the
  local clock time for a distributed agent~$i \in \until{N}$ as
  \begin{align*}
    s_i(t) = a_i t + b_i,
  \end{align*}
  where~$a_i > 0$ and~$b_i \in \real$ represent the unknown drift and
  bias of clock~$i$, respectively. In order to synchronize the clocks
  with varying and unknown drifts and biases, we define a virtual
  clock
  \begin{align*}
    T_i(t) = \alpha_i(s_i(t)) s_i(t) + \beta_i(s_i(t)),
  \end{align*}
  where~$\alpha_i$ and~$\beta_i$ are the controlled drift and bias of
  node~$i$, respectively, whose dynamics are to be designed. The goal
  is then to perform consensus on the virtual clock variables~$T_i$
  such that $\TwoNorm{T_i - T_j} \rightarrow 0$ for all pairs of
  agents by constantly sharing these values as they update their
  controlled drift and bias.  } Instead, recent works have studied
applying event-triggered coordination to these problems to reduce the
communication required by the agents to synchronize their
clocks~\citep{YK-HI:15,ZC-DL-YH-CT:15,EG-SM-YC-DWC:17}.

\paragraph*{Distributed optimization.} 
{In many applications, it is of interest to solve optimization
  problems with separable objective functions  of the form
  \begin{align*}
    f(x) = \sum_{i=1}^N f_i(x,y_i),
  \end{align*}
  subject to various equality and/or inequality constraints, where~$x$
  is some variable or state, ~$y_i$ is some data or measurement local
  to agent~$i$, and $f_i$ is a function known to agent $i$,
  see~\citep{PW-MDL:09,AN:14-sv}. Regardless of the
  objective or constraints, a popular approach to solve this problem
  is by having agents constantly share information (either their local
  state $x_i$ or their about the global solution $x$) with their
  neighbors to cooperatively optimize the function while ensuring the
  global constraints are satisfied. Instead, several groups of
  researchers have studied the application of event-triggered
  coordination to these sorts of problems where the agents' decisions
  are coupled through the objective function, the constraints, or
  both.}  More specifically, agents employ event-triggered
communication to trade computation for reduced overall information
exchange in finding a solution of the optimization
problem~\citep{MZ-CGC:10,SSK-JC-SM:15-auto,DR-JC:16,JL-WC:16,WC-WR:16}.

{
  \section{A Look Beyond}\label{se:beyond}

  We have mainly focused on providing an introduction and survey of
  event-triggered coordination strategies applied specifically to the
  multi-agent consensus problem, but this article can also be viewed
  as a broader tutorial on how to apply event-triggered coordination
  to networked systems in general by thinking of the consensus problem
  as a specific case study.  In fact, the work on event-triggered
  coordination for multi-agent consensus and networked systems in
  general has evolved in parallel over the last decade. For a tutorial
  introduction on event-triggered control of single plant systems we
  refer to~\citep{WPMHH-KHJ-PT:12}. Instead, we provide here a brief
  overview of distributed event-triggered control for networked
  systems in a more general context than just multi-agent
  consensus. Note that this is not meant to be a comprehensive survey
  of event-triggered coordination for networked systems in general but
  rather points to various other classes of networked problems that
  also benefit from ideas in event-triggered control; for this we
  refer the interested reader to~\citep{DT-SH:17}.  In this section,
  we provide a brief discussion on distributed event-triggered
  coordination of more general networked systems and finally take a
  look forward, identifying possible avenues for future research.

  \subsection{Distributed event-triggered control and stabilization}
  We consider here the distributed control and stabilization of
  interconnected systems.  These interconnections not only capture the
  ability to share information, but also might represent inherent
  dynamic coupling between the subsystems. The overall objective is
  the stabilization of the overall system through coordinated
  communication and control. Initial work formally introducing and
  addressing this problem was presented
  in~\citep{XW-MDL:08,XW-MDL:08b}, where each subsystem's dynamics are
  described by
  \begin{align}
    \dot{x}_i (t)= A_ix_i(t) + B_iu_i(t) + \sum_{j\in\mathcal{N}_i}
    A_{ij}x_j(t) , \label{eq:DSlinear}
  \end{align}
  for $i = \until{N}$, where $x_i\in \mathbb{R}^{n_i}$ and $u_i\in
  \mathbb{R}^{m_i}$ represent the state and the control input of
  subsystem $i$, respectively. The matrices $A_i\in \mathbb{R}^{n_i
    \times n_i}$, $B_i\in \mathbb{R}^{n_i \times m_i}$, $A_{ij}\in
  \mathbb{R}^{n_i \times n_j}$, represent the state, input, and
  coupling matrices, respectively. The starting point is the
  availability of ideal controllers
  \begin{align*}
    u_i(t) = K_i x_i + \sum_{j \in \NN_i} L_{ij} x_j,
  \end{align*}
  where the decoupling gains~$L_{ij}$ can be chosen such that~$B_i
  L_{ij} = -A_{ij}$ and~$K_i$ has been designed such~$x = 0$ is
  asymptotically stable. However, since implementing this solution
  requires each subsystem to have continuous access to the state of
  neighboring subsystems, the goal now is to instead design a
  distributed event-triggered coordination strategy to solve the
  problem.

  Similar to our setup for consensus, each subsystem must now
  determine for itself when to broadcast its state to neighboring
  subsystems. Consequently, the actual control input of agent~$i$ is
  given by
  \begin{align}\label{eq:stablecontrol}
    u_i(t) = K_i \hat{x}_i + \sum_{j \in \NN_i} L_{ij} \hat{x}_j,
  \end{align}
  where $ \hat{x}_i(t) = x_i(t_\ell^i) \text{ for } t \in
  [t_\ell^i,t_{\ell+1}^i)$ is the last broadcast state of
  subsystem~$i$ at any given time~$t \in \realnonnegative$. Just as in
  the consensus case, we define~$e_i(t)=\hat{x}_i(t) - x_i(t)$ as the
  error between the last broadcast state of subsystem~$i$ and its
  current state. Following our discussion in Section~\ref{se:primer},
  we are interested in designing a distributed event-triggering
  condition~$f_i(\cdot)$ such that the closed-loop dynamics
  of~\eqref{eq:DSlinear} with control input~\eqref{eq:stablecontrol}
  guarantees that the system still asymptotically converges.

  To achieve this, under suitable conditions,~\cite{XW-MDL:08} propose
  the trigger
  \begin{align*}
    f_i(e_i,x_i) \triangleq \beta_i \TwoNorm{e_i}^2 - \alpha_i
    \TwoNorm{x_i}^2 \leq 0 ,
  \end{align*}
  for some positive constants~$\beta_i$ and~$\alpha_i$.  This work was
  extended in~\citep{XW-MDL:11,CDP-RS-FW:13} to nonlinear systems of
  the form
  \begin{align*}
    \dot{x}_i (t)= F_i( x_i, \{x_j\}_{j \in \NN_i}, u_i).
  \end{align*}
  Later works~\citep{MG-DL-JS-SD-KHJ:12,MG-DVD-KHJ-JS-SD:13} have
  further extended these results to consider robustness issues such as
  the presence of network delays and packet dropouts in the networked
  stabilization problem.

  \subsection{Decentralized event-triggered control over wireless
    sensor/actuator networks}
  A similar problem related to distributed event-triggered control was
  presented in~\citep{MMJ-PT:10} and \citep{RP-PT-DN-AA:11}, where the
  states of a nonlinear system
  \begin{align*}%  \label{eq:DSnonlinearOR}
    \dot{x}(t)=F(x(t),u(t))
  \end{align*}
  are measured by individual, decentralized sensors. More
  specifically, different components of the state of a single system
  is sampled by different sensors. As usual, the starting point here
  is the assumption that a stabilizing controller~$u(t) = K(x(t))$
  exists and is known. The actually used control input is then of the
  form $u(t)=K(\hat{x}(t))$, where~$\hat{x}(t)$ represents the most
  up-to-date information about the true state~$x(t)$. The goal is then
  for any given sensor to determine conditions for sharing its
  information based only on its local measurements in order to
  guarantee stabilization of the overall system.

  \subsection{Distributed event-triggered control for output
    feedback systems}

  Throughout this article we have assumed that individual agents or
  subsystems have access to their own exact state. Instead, many
  works~\citep{PT-NC:12,MCFD-WPMHH:12,WPMHH-MCFD-ART:13,PT-NC:14,FF-SG-DN-LZ:14,MMJ-MC:14,HY-PA:12,DT-RF:13,HZ-GF-HY-QC:14b,YC-MF-DD:16,WH-LL:17,HY-PJA:14,WH-LL-GF:16,TL-MC-DJH:14,BZ-XL-TH-HL-GC:16,VSD-DPB-WPMHH:17,WH-LL-GF:17,MA-RP-JD-DN:17,MSM-MS-ME:16}
  consider the case of systems with output feedback, e.g.,
  \begin{align*}
    \dot{x}_i(t)&=A_ix_i(t) + B_iu_i (t) ,
    \\
    y_i(t)&=C_ix_i(t) ,
  \end{align*}
  where each element of the output of the system is sampled by a
  different sensor. Depending on the specific application, the goals
  are similar to the problems we have already discussed; e.g.,
  synchronizing the outputs (rather than the states) or guaranteeing
  asymptotic convergence using only an output feedback controller
  (rather than state feedback). 
  
  Similar to the technical issues we come across regarding Zeno executions for the
  multi-agent consensus problem, these types of problems exhibit the exact same types
  of concerns due to the distributed and partial information available to the sensors~\cite{MCFD-WPMHH:12}.}

%\marginJC{No comment here on the Zeno aspect? I remember Reviewer 4
%  bringing up the issue of output feedback and Zeno, and how that is
%  connected to distributed systems (because an agent only measures an
%  output, rather than the whole state, of the network). It'd be nice
%  to gave a paragraph here discussing it.}

% {\color{green}
% \subsubsection{Stochastic event-triggered control}
% 
%% This entire article has only focused on deterministic triggering
%% strategies. However, it should be noted that stochastic
%% event-triggered control strategies may also be developed to solve
%% the problems considered in this article and elsewhere. More
%% specifically, it is possible that agents use probabilistic
%% triggering conditions rather than deterministic ones
%% 
% 
%\margin{I actually couldn't find a single work that does distributed
%  stochastic event-triggered control. Note this is different from
%  deterministic event-triggered strategies for stochastic control. In
%  either case I don't think we should mention this in the article...}
%%}
%
%\marginJC{How about then having a paragraph in future outlook about
%  this. Somebody should be looking at it! We could point to the
%  references brought up by Reviewer 4 as an initial point of entry
%  into the problem.}

{
\subsection{Future Outlook}\label{se:outlook}

This article has provided an introduction to the field of
event-triggered consensus. As consensus problems are widespread in
terms of networked systems, the exposition has brought up many of the
challenges and tools that are not specific to this particular problem,
but underlie network coordination tasks in general. In particular, we
have highlighted the importance of designing distributed
event-triggers that still guarantee global properties and the
technical difficulties that come with ensuring stability of the
resulting asynchronous executions. The focus on consensus has enabled
us to illustrate the motivation, the design methods, and the technical
challenges that arise in carrying this through. Here, we offer some
thoughts on additional lines of research that we believe are worthy of further
exploration in the future.

{
  \paragraph*{Stochastic event-triggers.} Through this article we have
  focused only on deterministic triggering strategies that determine
  precise times at which events should be triggered. Alternatively, a
  number of works~\citep{DH-YM-JW-SW-BS-LS:15,FDB-DA-FA:18,DA:13}
  consider the design of events that are triggered
  stochastically. While this area is still in under development, there
  seem to be some benefits that may be applicable to both the
  multi-agent consensus problem and networked systems more
  generally. For instance,
  \begin{enumerate}
  \item it may be easier to compare the performance of different
    algorithms against one another by considering average quantities
    or rates of transmissions rather than exact trajectories;
  \item they require less precise specifications, in that triggers can
    be more loosely defined since it is not critical that an event be
    triggered at an exact specified time; and
  \item it may be easier to ensure non-Zeno executions due to the less
    precise scheduling constraints.
  \end{enumerate}
  We are not yet aware of any works that study this for networked
  systems but believe this to be a worthwhile avenue for future
  exploration.  It should be noted that we view this differently from
  \emph{stochastic event-triggered control}, in which the idea of
  event-triggered control updates is applied to stochastic optimal
  controllers~\citep{YX-JPH:05,MR-JSB:07,AM-SH:09,AM-SH:10}.}

\paragraph*{Dynamic average consensus.} Throughout the article we have
generally assumed that the final convergence value of the entire
network depends on some function of the initial states of all
agents. However, we can also imagine extensions of these ideas to
{dynamic average consensus} problems, in which the group of agents is
expected to track some time-varying quantity about all the agents. For
instance, the agents may be trying to agree on the average temperature
in a room; but if the temperature in the room is changing, we need a
dynamic consensus algorithm to track the average temperate of the room
in real time.  Of course of the temperature is changing too quickly,
then we cannot expect accurate tracking of it in a distributed manner;
however, depending on assumptions on how quickly the quantity of
interest is changing, different bounds on the tracking error can be
provided using dynamic consensus
algorithms~\citep{RAF-PY-KML:06b,DPS-ROS-RMM:05-ifac,SSK-JC-SM:15-ijrnc}. In
fact, the dynamic average consenus problem is a natural generalization
of the static average consensus problem with applications in a wide
variety of areas.
 
Existing dynamic average consensus algorithms generally require
continuous or periodic communication between agents in order to
adequately track the quantity of interest.
Consequently,~\cite{SSK-JC-SM:15} considers the application of
event-triggered coordination to the dynamic average consensus problem
to relax this requirement. Similar to some algorithms presented in
this article,~\cite{SSK-JC-SM:15} presents event-triggered versions of
these algorithms that determine when communication should occur
according to some time- or state-varying thresholds over multi-agent
networks with single integrator and dynamic topologies. Convergence is
only guaranteed to within a neighborhood of the time-varying average
signal.
% control update and information push, continuous event detection,
% and two types of trigger dependence: constant threshold and
% state-dependent.

\paragraph*{Cloud-based event-triggered control.}

Throughout this article we have assumed various forms of peer-to-peer
communication; in all cases the messages were directly sent from one
agent to another. Instead, some recent works consider scenarios where
agents communicate indirectly through the use of the
cloud~\citep{AA-DL-DVD-KHJ:15,SLB-CN-GJP:16,AA-DL-DVD-KHJ:17}.
{More specifically, we have thus far assumed that when an agent
  decides to send information to another, it is \emph{passively}
  received by the receiving agent. It is possible that the packet can
  be dropped or delayed, but in either case the receiving agent does
  not have an active part in receiving the message. Instead, using a
  cloud-based communication model, an agent~$i$ can only publish
  things intended for other agents to a cloud repository, rather than
  sending it directly to them. In this case, until a receiving
  agent~$j$ \emph{actively} decides to connect to the cloud and
  download the information that is available, it will not even be
  aware of a pending message from the transmitting agent~$i$. And in
  turn, when the message is received, it reveals information about the
  past state and plans of agent $i$, and likely only partial
  information about the present ones. This raises a plethora of
  interesting opportunities for the design of promises among agents,
  and the technical analysis of the resulting coordination algorithms,
  see also~\citep{CN-JC:16-tac}. The use of the cloud also opens the
  possibility of network agents with limited capabilities taking
  advantage of high-performance computation capabilities to deal with
  complex dynamical processes.}

\paragraph*{Performance guarantees.}
Something conspicuously missing from the event-triggered consensus
literature in particular and the event-triggered control literature in
general are performance guarantees that quantify the benefit of this
approach over time-triggered or periodic implementations.  The
self-tuning nature of event-triggered control, where events are
\emph{tuned} to the execution of the task at hand, makes it appealing
at a conceptual and design level.  Periodic control, by contrast,
requires the a priori selection of stepsizes and the consideration of
worst cases in doing so.  Simulations have consistently shown the
promise of event-triggered algorithms over periodic ones in many
cases. Apart from some early work examining this
issue~\cite{KJA-BMB:99,KJA-BMB:02}, it is only recently that some
works have establishes results along these lines for systems with a
centralized controller or decision maker, see
e.g.,~\citep{DA-WPMHH:14,VSD-DPB-WPMHH:17,BAK-DJA-WPMHH:17,BAK-DJA-WPMHH:18,PO-JC:18}. There
is also some preliminary work to apply this to network
settings~\citep{CR-HS-KHJ:16,DPB-RG-WPMHHL:17,SHJH-DPB-WPMHH:17}, but
this area as a whole is still largely incomplete. We expect such
guarantees and characterizations of average communication rates to be
increasingly important as event-triggered coordination algorithms gain
further popularity.

Related to this is the need to have ways of comparing different
event-triggered algorithms that successfully achieve the same task in
order to understand which one is better and under what conditions.
Even in this article we have presented several algorithms that solve
the same problem, and even assume the same capabilities of the agents.
For example, Theorems~\ref{th:sampled-data1}
and~\ref{th:sampled-data2} both solve Problem~\ref{pr:event-control}
under the same assumptions of the agents' abilities and are able to
provide the same guarantee: asymptotic convergence (including non-Zeno
guarantees). However, this gives us no insight into the transient
performance of these algorithms. Consequently, it is unclear which
algorithm would be better to implement to solve a given problem.
Ultimately, the type of algorithm we want to design and implement
should be optimized for a certain task. For instance, in some cases it
may be desirable to reduce the communication burden of the network as
much as possible, but in other scenarios it may make more sense to use
frequent communication to yield faster convergence. In any case, we
are still in the need as a community for established metrics and
methods for comparing different algorithms against one another to be
able to optimize them to meet varying performance needs.

\section{Conclusions}\label{se:conclusions}

The application of event-triggered coordination to large-scale
networks is currently surging in interest due to the rising ubiquity
of interconnected cyber-physical systems. As the number of devices
connected to a shared network grows larger than previously dealt with
in the past, distributed time-triggered coordination strategies do not
scale well. Such limitations require a rethinking of the periodic
control paradigm towards opportunistic schemes, as the ones discussed
in this article, that take full advantage of the knowledge of the
agents, the environment, and the task to efficiently manage the
available resources.  

The work~\cite{KJA-BMB:02} concluded that
\begin{displayquote}
  ``There are an increasing number of applications where the
  assumption of constant sampling rate is no longer valid, typical
  examples are multi-rate sampling and networked systems. Lebesgue
  sampling (or event-triggered sampling) may be a useful
  alternative...

...It would be very attractive to have a system theory similar to the
one for periodic sampling.''
\end{displayquote}

We believe the field of event-triggered coordination is shaping up to
be precisely the type of theory needed in addressing the various
network problems that exist today, where aperiodic sampling,
communication, and control should be viewed as an opportunity rather
than a disturbance. We hope that researchers find this article useful
in developing a better understanding of event-triggered coordination
and designing other general networks that operate in an efficient and
adaptive fashion.

\section*{Acknowledgments}
The authors would like to thank the anonymous reviewers of the paper
for numerous suggestions that helped significantly improve the
presentation.
%% WE SHOULD ADD HERE FUNDING SOURCES
% This work was partially supported by AFOSR Award FA9550-15-1-0108.
  }

{
%{\scriptsize
%\bibliographystyle{plainnat}%%
%\bibliographystyle{abbrvnat}
\bibpunct{[}{]}{,}{n}{}{;}
%\bibliography{alias,ReferencesETCS_Cameron}
%%}

\appendix

\section*{Appendix}

\paragraph*{Proof of Theorem~\ref{th:main}}
Consider the Lyapunov function
\begin{align*}
  V(x) = \frac{1}{2} x^T L x .
\end{align*}
Then, given the dynamics~\eqref{eq:dynamics} and the continuous
control law~\eqref{eq:ideal},
\begin{align*}
  \dot{V}(x) = x^T L \dot{x} = -x^T L^T L x = - \TwoNorm{Lx}^2,
\end{align*}
where we have used the fact that~$L$ is symmetric. It is now clear
that using the continuous control law~\eqref{eq:ideal} we have
$\dot{V}(x) < 0$ for all $Lx \neq 0$. Using LaSalle's Invariance
Principle~\cite{HKK:02}, it can then be shown that
\begin{align*}
  x(t) \rightarrow \{ Lx = 0 \} = \{ x_i = x_j \hspace*{1ex} \forall i,j \in \until{N} \}
\end{align*}
as $t \rightarrow \infty$. Combining this with the fact that the sum
of all states is an invariant quantity concludes the proof,
\begin{align*}
  \frac{d}{dt} \left( \ones{N}^T x(t) \right) = \ones{N}^T \dot{x}(t)
  = -\ones{N}^T L x(t) = 0 .
\end{align*}
\hfill $\blacksquare$

\paragraph*{Proof of Theorem~\ref{th:broadcast1}}
Consider again the Lyapunov function
\begin{align*}
V(x) = \frac{1}{2} x^T L x. 
\end{align*}
Then, we see that the trigger~\eqref{eq:decentraltrigger2} guarantees
that~\eqref{eq:enforcethis} is satisfied at all times. Combined with
the closed-loop dynamics~\eqref{eq:decentraldynamics2}, we have
\begin{align*}
\dot{V} \leq \sum_{i=1}^N (\sigma_i - 1) (1 - a |\NN_i|)\hat{z}_i^2,
\end{align*}
which is strictly negative for all~$L\hat{x} \neq 0$. Similar to the
conclusion in the proof of Theorem~\ref{th:main}, we can similarly
show that $\hat{x} \rightarrow \{ L \hat{x} = 0 \}$. By noticing that
$\hat{x}_i$ is simply a sampled subset of the trajectory of~$x_i$, we
have that $x \rightarrow \{ L x = 0 \}$.  Finally, combining this
again with the fact that the sum of all states is an invariant
quantity concludes the proof.

\hfill $\blacksquare$

\paragraph*{Proof of Lemma~\ref{le:bound}}
Given the Lyapunov function
\begin{align*}
V(x) = \frac{1}{2} (x - \bar{x} \mathbf{1})^T (x - \bar{x} \mathbf{1})
\end{align*}
and the closed-loop dynamics~\eqref{eq:decentraldynamics2}, we have
\begin{align*}
\dot{V} = x^T \dot{x} - \bar{x}^T \dot{x} = - x^T L
  \hat{x} - \bar{x} L \hat{x} = - x^T L \hat{x}. 
\end{align*}
Recalling $e_i(t) = \hat{x}_i(t) - x_i(t)$, we can expand this out to
\begin{align*}
  \dot{V} &= -\hat{x}^T L \hat{x} + e^T L \hat{x}
  \\
  &= -\sum_{i=1}^N \sum_{j \in \NN_i} \left( \frac{1}{2} (\hat{x}_i -
    \hat{x}_j)^2 - e_i (\hat{x}_i-\hat{x}_j) \right) .
\end{align*}

Using Young's inequality for each product we can bound (see~\cite{CN-JC:16-auto}
for why this choice)
\begin{align*}
  e_i(\hat{x}_i - \hat{x}_j) \leq e_i^2 + \frac{1}{4} (\hat{x}_i -
  \hat{x}_j)^2
\end{align*}
which yields
\begin{align*}
  \dot{V} &\leq -\sum_{i=1}^N \sum_{j \in \NN_i} \left( \frac{1}{2}
    (\hat{x}_i - \hat{x}_j)^2 - e_i^2 - \frac{1}{4} (\hat{x}_i -
    \hat{x}_j)^2 \right)
  \\
  &= - \sum_{i=1}^N \sum_{j \in \NN_i} \left(
    \frac{1}{4}(\hat{x}_i-\hat{x}_j)^2 - e_i^2 \right)
  \\
  &= \sum_{i=1}^N e_i^2 | \NN_i | - \sum_{j \in \NN_i} \left(
    \frac{1}{4}(\hat{x}_i-\hat{x}_j)^2 \right) .
\end{align*}
\hfill $\blacksquare$

\paragraph*{\textbf{Proof of Theorem~\ref{th:time}}}
Let $\delta(t) = x(t) - \bar{x} \mathbf{1}$, where $\bar{x} = \frac{1}{N}
\sum_{i=1}^N x_i(0)$ is the average of all initial conditions. Then,
$\dot{\delta}(t) = - L \delta(t) - L e(t)$, yielding
\begin{align*}
  \delta(t) = e^{-Lt} \delta(0) - \int_0^t e^{-L(t-s)} L e(s) ds .
\end{align*}
Taking norms,
\begin{align*}
  \TwoNorm{\delta(t)} &\leq \TwoNorm{ \delta(0) e^{-Lt}} + \int_0^t
  \TwoNorm{ e^{-L(t-s)} Le(s) } ds
  \\
  &\leq e^{- \lambda_2(L) t} \TwoNorm{\delta(0)} + \int_0^t
  e^{-\lambda_2(L)(t-s)} \TwoNorm{L e(s)} ds, 
\end{align*}
where the second inequality follows from~\cite[Lemma 2.1]{GSS-DVD-KHJ:13}.

Using the condition
\begin{align*}
  | e_i(t) | \leq c_0 + c_1 e^{-\alpha t},
\end{align*}
it follows that
\begin{align*}
  \TwoNorm{\delta(t)} &\leq e^{-\lambda_2 t} \TwoNorm{\delta(0)} +
  \TwoNorm{L} \sqrt{N} \int_0^t e^{-\lambda_2(t-s)} (c_0 + c_1
  e^{-\alpha s}) ds
  \\
  &= e^{-\lambda_2 t} \left( \TwoNorm{\delta(0)} - \TwoNorm{L}
    \sqrt{N} \left( \frac{c_0}{\lambda_2} + \frac{c_1}{\lambda_2 -
        \alpha} \right) \right) \\ 
        &+ e^{-\alpha t} \frac{ \TwoNorm{L}
    \sqrt{N} c_1}{\lambda_2 - \alpha} + \frac{ \TwoNorm{L} \sqrt{N}
    c_0 }{\lambda_2}.
\end{align*}
The convergence result then follows by taking $t \rightarrow \infty$. 
\hfill $\blacksquare$

\paragraph*{\textbf{Proof of Theorem~\ref{th:sampled-data2}}}
Consider the Lyapunov function
\begin{align*}
  V(x) = \frac{1}{2}(x - \bar{x} \mathbf{1})^T ( x - \bar{x}
  \mathbf{1}).
\end{align*}
Following the discussing after Lemma~\ref{le:bound}, we know that
when~\eqref{eq:enforcethisalways} is satisfied, we have~$\dot{V}$ is
strictly negative for all~$L \hat{x} \neq 0$. However, since the
agents can now only evaluate the trigger~\eqref{eq:trigger3} at the
sampling times under the periodic event-triggered coordination
algorithm presented in Theorem~\ref{th:sampled-data2}, we lose the
guarantee that~$\dot{V} \leq 0$ at all times. Thus, we must now
analyze what happens to the Lyapunov function~$V$ in between these
sampling times. Explicitly considering $t \in [t_{\ell'},
t_{\ell'+1})$, note that
\begin{align*}
  e(t) = e(t_{\ell'}) + (t-t_{\ell'}) L \hat{x}(t_{\ell'}) .
\end{align*}
Substituting this expression into $\dot{V}(t) = -\hat{x}^T(t) L
\hat{x}(t) + e^T(t) L \hat{x}(t)$, we obtain
\begin{multline*}
  \dot{V}(t) = -\hat{x}^T(t_{\ell'}) L \hat{x}(t_{\ell'}) + e^T(t_{\ell'}) L
  \hat{x}(t_{\ell'})
  \\
  + (t - t_{\ell'}) \hat{x}^T(t_{\ell'}) L^T L \hat{x}(t_{\ell'}) ,
\end{multline*}
for all $t \in [t_{\ell'}, t_{\ell'+1})$.  For a simpler exposition,
we drop all arguments referring to time $t_{\ell'}$ in the sequel.
Then, using~\eqref{eq:enforcethisalways} to bound~$e(t_\ell')$, we
can show
\begin{align*}
  \dot{V}(t) \leq \sum_{i=1}^N \frac{\sigma_i - 1}{4} \sum_{j \in
    \NN_i} (\hat{x}_i - \hat{x}_j)^2 + (t-t_{\ell'})
  \hat{x}^T L^T L \hat{x} .
\end{align*}
Note that the first term is exactly what we have when we are able to
monitor the trigger continuously~\eqref{eq:thenthis}.

Using the fact that $\left( \sum_{k=1}^p y_k \right)^2 \leq p
\sum_{k=1}^p y_k^2$ (which follows directly from the Cauchy-Schwarz
inequality), we bound
\begin{align*}
  \hat{x}^T L^T L \hat{x} & = \sum_{i=1}^N \left( \sum_{j \in \NN_i}
   (\hat{x}_i-\hat{x}_j) \right)^2 
  \\
  & \le \sum_{i=1}^N |\NN_i| \sum_{j \in \NN_i } (\hat{x}_i -
  \hat{x}_j)^2  .
%  \\
%  & = |\NN_i | \sum_{i=1}^N \sum_{j \in \NN_i} (\hat{x}_i -
%  \hat{x}_j)^2 .
\end{align*}
%where $w^{\max}_i = \max_{j \in \NN_i^\text{out}} w_{ij}$. 
% Using
% \begin{align}
%   \hat{x}^T L^T L \hat{x} \leq C \hat{x}^T L \hat{x} ,
% \end{align}
% where
% \begin{align*}
%   C = \frac{\lambda_\text{max}(L^T L)}{\lambda_2( L_s )},
% \end{align*}
Hence, for $t \in [t_{\ell'}, t_{\ell'+1})$,
\begin{align*}
  \dot{V}(t) & \leq \sum_{i=1}^N \Big( \frac{\sigma_i - 1}{4} + h
  |\NN_\text{max} |^2 \Big) \sum_{j \in \NN_i} (\hat{x}_i -
  \hat{x}_j)^2
  \\
  & \leq \Big( \frac{\sigma_{\max}-1}{2} + 2 h
   |\NN_\text{max} |^2 \Big) \hat{x}^T L \hat{x} .
\end{align*}
Then, by using~\eqref{eq:bound2}, it can be shown that there exists $\BB > 0$
such that
\begin{align*}
  \dot{V}(t) & \leq \frac{1}{2 \BB }\Big( \sigma_{\max} + 4 h 
  |\NN_\text{max} |^2 - 1 \Big) V(x(t)),
\end{align*}
which implies the result. See~\cite{CN-JC:16-auto} for more details.
\hfill $\blacksquare$

\paragraph*{\textbf{Proof of Theorem~\ref{th:state-directed}}}

Consider the Lyapunov function
\begin{align*}
V(x) = \frac{1}{2}x^T L x .
\end{align*}
Then, we see that the trigger~\eqref{eq:trigger3} ensures that~\eqref{eq:enforce3} is satisfied at all times. Then, leveraging Lemma~\ref{lemma:bound}, we have
\begin{align*}
  \dot{V} \leq \sum_{i=1}^N \frac{\sigma_i - 1}{4} \sum_{j \in
    \NN_i^\text{out}} w_{ij} (\hat{x}_i - \hat{x}_j)^2,
\end{align*}
which is strictly negative for all~$L\hat{x} \neq 0$. Following the
discussion in the proof of Theorem~\ref{th:broadcast1}, we have that
$x \rightarrow \{Lx = 0\}$. Finally, because the graph is
weight-balanced the sum of all states is an
invariant quantity,
\begin{align*}
  \frac{d}{dt} \left( \mathbf{1}^T_N x(t) \right) = \mathbf{1}_N^T
  \dot{x}(t) = -\mathbf{1}_N^T L \hat{x}(t) = 0,
\end{align*}
which concludes the proof. \hfill $\blacksquare$

\paragraph*{\textbf{Proof of Theorem~\ref{th:DI-time}}}
Let $\delta(t) = x(t) - \bar{x} \mathbf{1}$, where $\bar{x} = \frac{1}{N}
\sum_{i=1}^N x_i(0)$ is the average of all initial conditions. Then,
$\dot{\delta}(t) = - L \delta(t) - L e(t)$, yielding
\begin{align*}
  \delta(t) = e^{-Lt} \delta(0) - \int_0^t e^{-L(t-s)} L e(s) ds .
\end{align*}
Taking norms,
\begin{align*}
  \TwoNorm{\delta(t)} &\leq \TwoNorm{ \delta(0) e^{-Lt}} + \int_0^t
  \TwoNorm{ e^{-L(t-s)} Le(s) } ds
  \\
  &\leq e^{- \lambda_2(L) t} \TwoNorm{\delta(0)} + \int_0^t
  e^{-\lambda_2(L)(t-s)} \TwoNorm{L e(s)} ds, 
\end{align*}
where the second inequality follows from~\cite[Lemma 2.1]{GSS-DVD-KHJ:13}.

Using the condition
\begin{align*}
  | e_i(t) | \leq c_0 + c_1 e^{-\alpha t},
\end{align*}
it follows that
\begin{align*}
  \TwoNorm{\delta(t)} &\leq e^{-\lambda_2 t} \TwoNorm{\delta(0)} +
  \TwoNorm{L} \sqrt{N} \int_0^t e^{-\lambda_2(t-s)} (c_0 + c_1
  e^{-\alpha s}) ds
  \\
  &= e^{-\lambda_2 t} \left( \TwoNorm{\delta(0)} - \TwoNorm{L}
    \sqrt{N} \left( \frac{c_0}{\lambda_2} + \frac{c_1}{\lambda_2 -
        \alpha} \right) \right) \\
        & \hspace*{2ex} + e^{-\alpha t} \frac{ \TwoNorm{L}
    \sqrt{N} c_1}{\lambda_2 - \alpha} + \frac{ \TwoNorm{L} \sqrt{N}
    c_0 }{\lambda_2}.
\end{align*}
The convergence result then follows by taking $t \rightarrow \infty$. 
See~\citep{GSS-DVD-KHJ:13} for details on excluding Zeno behavior.
\hfill $\blacksquare$

\paragraph*{Proof of Lemma~\ref{th:NS_LD}}
Define $x=[x_1, ..., x_N]^T$. Then, using the Kronecker product, the
dynamics of the overall system can be expressed as follows
\begin{align*}
  %\label{eq:OAsys}
  \dot{x} = (\bar{A}+\bar{B})x
\end{align*}
where $\bar{A}=I_N\otimes A$ and $\bar{B}=c{L}\otimes BF$. There
exists a similarity transformation $S$ such that ${L}_J=S^{-1}{L}S$ is
in Jordan canonical form. Define $\bar{S}=S\otimes I_n$ and calculate
the following:
\begin{align*}
      \bar{S}^{-1}(\bar{A}-\bar{B})\bar{S}  &= \bar{S}^{-1}\bar{A}\bar{S} + \bar{S}^{-1} ({L}\otimes BF) \bar{S}   \nonumber   \\
      & = I_N\otimes A+ c {L}_J \otimes BF.    %\label{eq:LDst}
\end{align*}
By applying the similarity transformation we obtain that the eigenvalues of $\bar{A}+\bar{B}$ are given by the eigenvalues of $A+c\lambda_j BF$, where $\lambda_j=\lambda_j({L})$. \hfill $\blacksquare$

\end{document}